\documentstyle[draft,leqno,epsf,version]{aism_ftsm2}
\makeatletter
\renewcommand{\theequation}{\thesection\arabic{equation}}
\@addtoreset{equation}{section}
\makeatother

\newcommand{\bdm}{\begin{displaymath}}
\newcommand{\edm}{\end{displaymath}}
\newcommand{\eqref}[1]{(\ref{#1})}
\renewcommand{\(}{\left(}
\renewcommand{\)}{\right)}
\newcommand{\Cost}{{\cal C}}
\newcommand{\train}{{\cal T}}
\newcommand{\Dict}{{\cal D}}
\newcommand{\ndim}{n}
\newcommand{\ndimlev}{{n_0}}
\newcommand{\maxlev}{K}
\newcommand{\nsamp}{N}
\newcommand{\LF}{L}
\newcommand{\N}{{\bf N}}

\newcommand{\R}{{\bf R}}
\newcommand{\ba}{{\mbox{\boldmath $a$}}}
\newcommand{\bb}{{\mbox{\boldmath $b$}}}
\newcommand{\be}{{\mbox{\boldmath $e$}}}

\newcommand{\bw}{{\mbox{\boldmath $w$}}}
\newcommand{\bx}{{\mbox{\boldmath $x$}}}
\newcommand{\by}{{\mbox{\boldmath $y$}}}
\newcommand{\bX}{{\mbox{\boldmath $X$}}}
\newcommand{\bY}{{\mbox{\boldmath $Y$}}}
\newcommand{\bzero}{{\mbox{\boldmath $0$}}}
\newcommand{\bone}{{\mbox{\boldmath $1$}}}
\newcommand{\define}{\stackrel{\Delta}{=}}
\newcommand{\dd}[1]{\, {\rm d}{#1}}

\newcommand{\OG}{{\rm O}}
\newcommand{\SO}{{\rm SO}}
\newcommand{\SL}{{\rm SL}}
\newcommand{\GL}{{\rm GL}}
\newcommand{\e}{{\rm e}}
\newcommand{\Prob}{{\rm Pr}}
\newcommand{\subsp}{\Omega}
\excludeversion{hide}

\volumenumber{54} \vnumber{1} \publishyear{2002}     
 \frompage{1} \topage{4}        
\shorttitle{Sparsity vs. Independence}        
\shortauthor{B. B\'enichou and N. Saito} 
\title{Sparsity vs. Statistical Independence in Adaptive Signal Representations: A Case Study of the Spike Process\footnote{This research was partially supported by NSF DMS-99-73032, DMS-99-78321, and ONR YIP N00014-00-1-046.}}
\author{\sc 
Bertrand B\'enichou$^1$ and Naoki Saito$^2$
}    
\affiliation{
$^1$Ecole Nationale Sup\'erieure des T\'el\'ecommunications, 46, rue Barrault, 75634 Paris cedex 13 France\\
$^2$Department of Mathematics,          
University of California, Davis, CA 95616 USA\\
}
\abstract{                                             
Finding a basis/coordinate system that can efficiently represent 
an input data stream by viewing them as realizations of a stochastic process
is of tremendous importance in many fields including data compression and
computational neuroscience.
Two popular measures of such efficiency of a basis are 
sparsity (measured by the expected $\ell^p$ norm, $0 < p \leq 1$)
and statistical independence (measured by the mutual information).
Gaining deeper understanding of their intricate relationship, however, remains
elusive.
Therefore, we chose to study a simple synthetic stochastic process called 
the spike process, which puts a unit impulse at a random location in 
an $\ndim$-dimensional vector for each realization.
For this process, we obtained the following results:
1) The standard basis is the best both in terms of sparsity and statistical
independence if $\ndim \geq 5$ and the search of basis is restricted within
all possible orthonormal bases in $\R^\ndim$;
2) If we extend our basis search in all possible invertible linear 
transformations in $\R^\ndim$, then the best basis in statistical
independence differs from the one in sparsity;
3) In either of the above, the best basis in statistical independence is
not unique, and there even exist those which make the inputs completely dense;
4) There is no linear invertible transformation that achieves the true
statistical independence for $\ndim > 2$.
}
\begin{document}

\maketitle

\section{Introduction}
\label{sec:intro}
What is a good coordinate system/basis to efficiently represent a given set of
images?  We view images as realizations of a certain complicated 
stochastic process whose probability density function (pdf)
is not known a priori.
\emph{Sparsity} is important here since this is a measure of how well one can 
compress the data. A coordinate system producing a few large coefficients and 
many small coefficients has high sparsity for that data.
The sparsity of images relative to a coordinate system is often measured by
the expected $\ell^p$ norm of the coefficients where $0 < p \leq 1$.
\emph{Statistical independence} is also important since 
statistically independent coordinates do not interfere with each other
(no crosstalk, no error propagation among them).
The amount of statistical dependence of input images relative to a
coordinate system is often measured by the so-called mutual information, 
which is a statistical distance between the true pdf and the product of 
the one-dimensional marginal pdf's.

Neuroscientists have become interested in efficient representations
of images, in particular, images of natural scenes such as trees, rivers,
mountains, etc., since our visual system effortlessly reduces the amount of
visual input data without losing the essential information contained in them.
Therefore, if we can find what type of basis functions
are sparsifying the input images or are providing us with the statistically
independent representation of the inputs, then that may shed light on the
mechanisms of our visual system.  Olshausen and Field (1996, 1997)
pioneered such studies using computational 
experiments emphasizing the sparsity.  Immediately after their experiments, 
Bell and Sejnowski (1997), van Hateren and van der Schaaf (1998)
conducted similar studies using 
the statistical independence criterion.
Surprisingly, these results suggest that both sparsity and independence 
criteria tend to produce oriented Gabor-like functions, which are similar 
to the receptive field profiles of the neurons in our primary visual cortex.
However, the relationship between these two criteria has not been 
understood completely.

These experiments and observations inspired our study in this paper.
We wish to deepen our understanding of this intricate relationship.
Our goal here, however, is more modest in that we only study the so-called 
``spike'' process, a simple synthetic stochastic
process, which puts a unit impulse at a random location in an 
$\ndim$-dimensional vector for each realization.
It is important to use a simple stochastic process first since 
we can gain insights and make precise statements in terms of theorems.
By these theorems, we now understand what are the precise conditions for
the sparsity and statistical independence criteria to select the same basis
for the spike process.  In fact, we prove the following facts. 
\begin{itemize}
\item The standard basis is the best both in terms of sparsity and statistical
independence if $\ndim \geq 5$ and the search of a basis is restricted
within all possible orthonormal bases in $\R^\ndim$.
\item If we extend our basis search in all possible invertible
linear transformations in $\R^\ndim$, then the best basis in statistical
independence differs from the standard basis, which is the best in sparsity.
\item In either of the above, the best basis in statistical independence 
is not unique, and there even exist those which make the inputs completely 
dense;
\item There is no linear invertible transformation that achieves the true
statistical independence for $\ndim > 2$.
\end{itemize}
These results and observations hopefully lead to deeper understanding of 
the efficient representations of more complicated stochastic processes such as
natural scene images.

More information about other stochastic processes, such as
the ``ramp'' process (another simple yet important stochastic
process), can be found in Saito et al.\ (2000, 2001), which also contain our
numerical experiments on natural scene images.

The organization of this paper is as follows.
In Section~2, we set our notations and terminology.
Then in Section~3, we precisely define how to quantitatively measure 
the sparsity and statistical dependence of a stochastic process relative to 
a given basis.  Using a very simple example, Section~4 demonstrates that 
the sparsity and statistical independence are two clearly different
concepts.  Section~5 presents our main results.
We prove these theorems in Section~6 and Appendices.
Finally, we discuss the implications and further directions in Section~7.

\section{Notations and Terminology}
\label{sec:notation}
Let us first set our notation and the terminology of basis
dictionaries and best bases.
Let $\bX \in \R^\ndim$ be a random vector with some unknown pdf $f_\bX$.
Let us assume that the available data $\train=\{\bx_1,\ldots,\bx_\nsamp\}$
were independently generated from this probability model.  The set
$\train$ is often called the training dataset.
Let $B = \( \bw_1, \ldots, \bw_\ndim \) \in \OG(\ndim)$ (the group of
orthonormal transformations in $\R^\ndim$) or $\SL^\pm(\ndim,\R)$ (the group
of invertible volume-preserving transformations in $\R^\ndim$, i.e.,
their determinants are $\pm 1$).
The best-basis paradigm, Coifman and Wickerhauser (1992), Wickerhauser (1994),
Saito (2000), is to find a basis $B$ or a subset of basis 
vectors such that the features (expansion coefficients) $\bY = B^{-1} \bX$
are useful for the problem at hand (e.g., compression, modeling, 
discrimination, regression, segmentation) in a computationally fast manner.
Let $\Cost(B \, | \, \train)$ be a numerical measure of \emph{deficiency} 
or \emph{cost} of the basis $B$ given the training dataset $\train$ for 
the problem at hand.
For very high-dimensional problems, we often restrict our search within
the basis dictionary $\Dict \subset
\SL^\pm(\ndim,\R)$, such as the orthonormal or biorthogonal wavelet packet 
dictionaries or local cosine or Fourier dictionaries where we never need to
compute the full matrix-vector product or the matrix inverse for analysis and 
synthesis.
Under this setting,
$B_\star = \arg \min_{B \in \Dict} \Cost(B \, | \, \train)$ is called
the \emph{best basis} relative to the cost $\Cost$ and the training dataset
$\train$.  Section~\ref{sec:review-Haar-Walsh} reviews the concept of 
the basis dictionary and the best-basis algorithm in details.

We also note that $\log$ in this paper implies $\log_2$, unless stated
otherwise.

\section{Sparsity vs. Statistical Independence}
\label{sec:spar.vs.indep}
The concept of sparsity and that of statistical independence are
intrinsically different.  Sparsity emphasizes the issue of compression
directly, whereas statistical independence concerns the relationship among 
the coordinates.
Yet, for certain stochastic processes, these two are intimately related, and
often confusing.  For example, Olshausen and Field (1996, 1997) emphasized
the sparsity as the basis selection 
criterion, but they also assumed the statistical independence of the 
coordinates.
Bell and Sejnowski (1997) used the statistical independence criterion and 
obtained the basis functions similar to those of Olshausen and Field.
They claimed that they did not impose the sparsity explicitly and such sparsity
\emph{emerged} by minimizing the statistical dependence among the coordinates.
These motivated us to study these two criteria.

First let us define the measure of sparsity and that of statistical 
independence in our context.
\subsection{Sparsity}
Sparsity is a key property as a good coordinate system for compression.
The true sparsity measure for a given vector $\bx \in \R^\ndim$ is 
the so-called $\ell^0$ quasi-norm
which is defined as
\bdm
\| \bx \|_0 \define \#\{ i \in [1,\ndim] : x_i \neq 0 \},
\edm
i.e., the number of nonzero components in $\bx$.
This measure is, however, very unstable for even small perturbation
of the components in a vector.  Therefore, a better measure is
the $\ell^p$ norm:
\bdm
\label{eq:ellp-norm}
\| \bx \|_p \define \( \sum_{i=1}^\ndim |x_i|^p \)^{1/p}, \quad 0 < p \leq 1.
\edm
In fact, this is a quasi-norm for $0 < p < 1$ 
since this does not satisfy the triangle
inequality, but only satisfies weaker conditions:
$\| \bx + \by \|_p \leq 2^{-1/p'} (\| \bx \|_p + \| \by \|_p)$ where
$p'$ is the conjugate exponent of $p$; and
$\| \bx + \by \|^p_p \leq \| \bx \|^p_p + \| \by \|^p_p$.
It is easy to show that $\lim_{p \, \downarrow \, 0} \| \bx \|_p^p = \| \bx \|_0$.
See Day (1940), Donoho (1994, 1998) for the details of the $\ell^p$ norm 
properties.

Thus, we can use the expected $\ell^p$ norm minimization as a criterion 
to find the best basis for a given stochastic process in terms of sparsity:
\begin{equation}
\label{eq:cp}
\Cost_{p}(B \, | \, \bX) = E \| B^{-1} \bX \|^p_p,
\end{equation}
The sample estimate of this cost given the training dataset $\train$ is
\begin{equation}
\label{eq:cp-samp}
\Cost_{p}(B \, | \, \train) = \frac{1}{\nsamp} \sum_{k=1}^\nsamp \| \by_k \|^p_p = \frac{1}{\nsamp} \sum_{k=1}^\nsamp \sum_{i=1}^\ndim | y_{i,k} |^p,
\end{equation}
where $\by_k=(y_{1,k},\ldots,y_{\ndim,k})^T=B^{-1}\bx_k$\, and $\bx_k$ is
the $k$th sample (or realization) in $\train$. 
We propose to use the minimization of this cost to select the \emph{best sparsifying basis} (BSB):
\bdm
\label{eq:min-cp}
	B_{p}= B_{p}(\train, \Dict) = \arg \min_{B \in \Dict} \Cost_{p}(B \, | \, \train).
\edm
\begin{aismrm}
It should be noted that \emph{the minimization of the $\ell^p$ norm can also be
achieved for each realization}. 
Without taking averages in \eqref{eq:cp-samp}, 
one can select the BSB $B_p=B_p(\{\bx_k\}, \Dict)$ for each 
realization $\bx_k \in \train$.  We can guarantee that
\bdm
\min_{B \in \Dict} \Cost_p (B \, | \, \{ \bx_k \}) \leq \min_{B \in \Dict} \Cost_p (B \, | \, \train) \leq \max_{B \in \Dict} \Cost_p (B \, | \, \{ \bx_k \}).
\edm
For highly variable or erratic stochastic processes, however, $B_p(\{\bx_k\}, \Dict)$
may significantly change for each $k$ and we need to store more information of
this set of $\nsamp$ bases if we want to use them to compress the entire
training dataset.  Whether we should adapt a basis per realization or
on the average is still an open issue.  See Saito et al.\ (2000, 2001)
for more details.
\end{aismrm}

\subsection{Statistical Independence}
The statistical independence of the coordinates of $\bY \in \R^\ndim$ means
\bdm
f_\bY(\by) = f_{Y_1}(y_1) f_{Y_2}(y_2) \cdots f_{Y_\ndim}(y_\ndim),
\edm
where $f_{Y_k}(y_k)$ is a one-dimensional marginal pdf.
The statistical independence
is a key property as a good coordinate system for compression and particularly
modeling because: 1) damage of one coordinate does not propagate to the others;
and 2) it allows us to model the $\ndim$-dimensional stochastic process of 
interest as a set of 1D processes.
Of course, in general, it is difficult to find a truly statistically 
independent coordinate system for a given stochastic process. 
Such a coordinate system may not even exist for a certain stochastic process.
Therefore, we should be satisfied with finding the least-statistically
dependent coordinate system within a basis dictionary.
Naturally, then, we need to measure the ``closeness'' of a coordinate
system $Y_1,\ldots,Y_\ndim$ to the statistical independence.
This can be measured by \emph{mutual information}
or relative entropy between the true pdf $f_\bY$ and the product of its
marginal pdf's:
\bdm
	I(\bY)  \define \int f_\bY(\by) \log \frac{f_\bY(\by)}{\prod_{i=1}^\ndim f_{Y_i}(y_i)} \dd{\by} \\
                =  - H(\bY) + \sum_{i=1}^\ndim H(Y_i),
\edm
where $H(\bY)$ and $H(Y_i)$ are the differential entropy of $\bY$ and $Y_i$
respectively:
\bdm
 H(\bY) = - \int f_\bY(\by) \, \log f_\bY(\by) \dd{\by}, \quad
 H(Y_i) = - \int f_{Y_i}(y_i) \, \log f_{Y_i}(y_i) \dd{y_i}.
\edm
We note that $I(\bY) \geq 0$, and $I(\bY) = 0$ if and only if
the components of $\bY$ are mutually independent. See Cover and Thomas (1991)
for more details of the mutual information.

Suppose $\bY = B^{-1} \bX$ and $B \in \GL(\ndim,\R)$ with $\det(B)=\pm 1$.
We denote such a group of matrices by $\SL^\pm(\ndim,\R)$.  Note that
the usual $\SL(\ndim,\R)$ is a subgroup of $\SL^\pm(\ndim,\R)$.
Then, we have
\bdm
I(\bY) = - H(\bY) + \sum_{i=1}^\ndim H(Y_i) = - H(\bX) + \sum_{i=1}^\ndim H(Y_i),
\edm
since the differential entropy is \emph{invariant}
under such a invertible volume-preserving linear transformation, i.e., 
\bdm
H(B^{-1}\bX)=H(\bX)+\log | \det(B^{-1})|= H(\bX),
\edm
because $| \det(B^{-1}) | = 1$.
Based on this fact, we proposed the minimization of the following cost 
function as the criterion to select
the so-called \emph{least statistically-dependent basis} (LSDB) 
in Saito (2001):
\begin{equation}
\label{eq:c-lsdb}
\Cost_{H}(B \, | \, \bX)  = \sum_{i=1}^\ndim H\((B^{-1}\bX)_i\)
= \sum_{i=1}^\ndim H(Y_i). 
\end{equation}
The sample estimate of this cost given the training dataset $\train$ is
\bdm
\label{eq:c-lsdb-samp}
\Cost_{H}(B \, | \, \train) = - \frac{1}{\nsamp}\sum_{k=1}^\nsamp \sum_{i=1}^\ndim \log \widehat{f}_{Y_i}(y_{i,k}),
\edm
where $\widehat{f}_{Y_i}(y_{i,k})$ is an empirical pdf of the coordinate
$Y_i$, which must be estimated by an algorithm such as the histogram-based
estimator with optimal bin-width search of Hall and Morton (1993).
Now, we can define the LSDB as
\begin{equation}
\label{eq:min-c-lsdb}
	B_{LSDB} = B_{LSDB}(\train, \Dict) = \arg \min_{B \in \Dict} \Cost_{H}(B \, | \, \train).
\end{equation}
We note that the differences between this strategy and the standard independent
component analysis (ICA) algorithms are: 1) restriction of the search
in the basis dictionary $\Dict$; and 2) approximation of the coordinate-wise
entropy.  For more details, we refer the reader to Saito (2001)
for the former and Cardoso (1999) for the latter.

Now we describe our analysis of some simple stochastic processes.

\section{Two-Dimensional Counterexample}
\label{sec:2Dunif}
\begin{figure}
\epsfxsize=0.7\textwidth
\centerline{\epsfbox{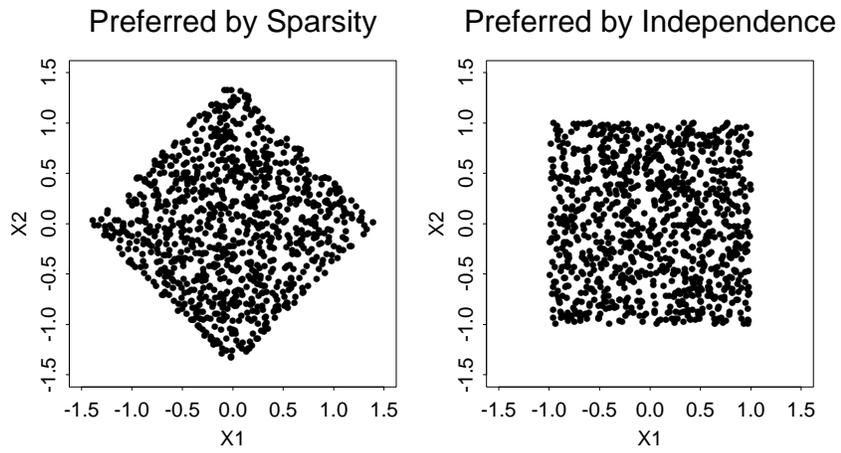}}
\caption{Sparsity and statistical independence prefer the different
coordinates.}
\label{fig:skew-noskew}
\end{figure}

This example clearly demonstrates the difference between the sparsity and
the statistical independence criteria.
Let us consider a simple process $\bX=(X_1,X_2)^T$ where $X_1$ and $X_2$ are
independently and identically distributed as the uniform random variable
on the interval $[-1,1]$.
Thus, the realizations of this process are distributed as the right-hand side
of Figure~\ref{fig:skew-noskew}.
Let us consider all possible rotations around the origin as
a basis dictionary, i.e., $\Dict=\SO(2,\R) \subset \OG(2)$.
Then, the sparsity and independence criteria select completely different
bases as shown in Figure~\ref{fig:skew-noskew}.
Note that the data points under the BSB coordinates (45 degree rotation)
concentrate more around the origin than the LSDB coordinates (with no rotation)
and this makes the data representation sparser.
This example clearly demonstrates that the BSB and the LSDB are different
in general.
One can also generalize this example to higher dimensions.

\section{The Spike Process}
\label{sec:spike}
\begin{figure}
\epsfxsize=0.7\textwidth
\centerline{\epsfbox{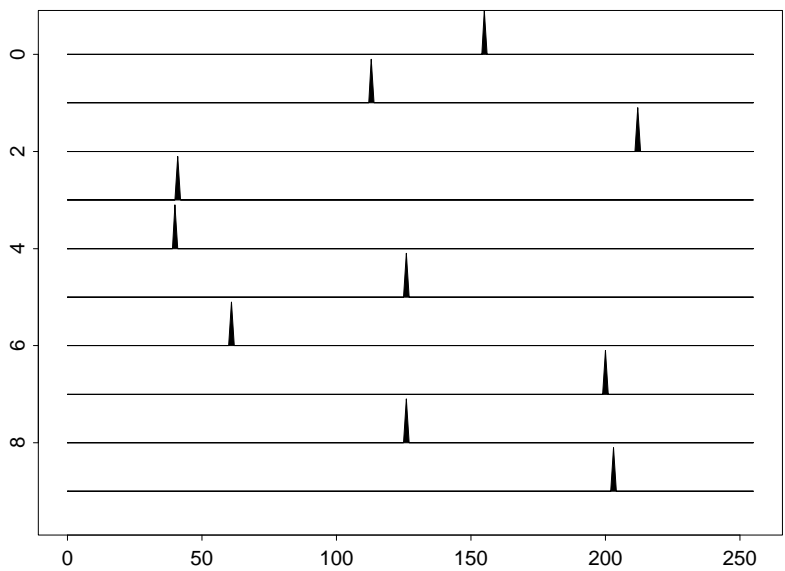}}
\caption{Ten realizations of the spike process ($\ndim=256$).}
\label{fig:spikes}
\end{figure}

An $\ndim$-dimensional \emph{spike process} simply generates the standard
basis vectors $\{ \be_j \}_{j=1}^\ndim \subset \R^\ndim$ in a random order,
where $\be_j$ has one at the $j$th entry and all the other entries are zero.
One can view this process as a unit impulse located at a
random position between $1$ and $\ndim$ as shown in Figure~\ref{fig:spikes}.

\subsection{The Karhunen-Lo\`eve Basis}
Let us first consider the Karhunen-Lo\`eve basis of this process from which
we can learn a few things.
\begin{aismpr}
\label{prop:spike-KLB}
The Karhunen-Lo\`eve basis for the spike process is any orthonormal basis 
in $\R^\ndim$ containing the ``DC'' vector $\bone_\ndim=(1,1,\ldots,1)^T$.
\end{aismpr}
This means that the KLB is not useful for this process.  This is because
the spike process is highly non-Gaussian.

\subsection{The Best Sparsifying Basis}
It is obvious that the standard basis is the BSB among $\OG(\ndim)$
by construction; an expansion of a realization of this process into any other
basis simply increases the number of nonzero coefficients.  
More precisely, we have the following proposition.
\begin{aismpr}
\label{prop:spike-sparse}
The BSB for the spike process is the standard basis
if $\Dict=\OG(\ndim)$ or $\SL^\pm(\ndim,\R)$.
If $\Dict=\GL(\ndim,\R)$, then it must be a scalar
multiple of the identity matrix, i.e., $aI_\ndim$ where $a$ is a nonzero
constant.
\end{aismpr}

\begin{aismrm}
\label{rem:permute-sign}
Note that when we say the basis is a matrix such as $aI_\ndim$,
we really mean that the column vectors of that matrix form the basis.
This also means that any permuted and/or sign-flipped (i.e., multiplied by
$-1$) versions of those column vectors also form the basis.  
Therefore, when we say the basis is a matrix $A$, we mean not only
$A$ but also its permuted and sign-flipped versions of $A$.
This remark also applies to all the propositions, lemmas, and
theorems below, unless stated otherwise.
\end{aismrm}

\subsection{Statistical Dependence and Entropy of the Spike Process}
\label{sec:true-entropy}
Before considering the LSDB of this process, let us note a few specifics
about the spike process.
First, although the standard basis is the BSB for this process, 
it clearly does not provide the statistically independent coordinates.
The existence of a single spike at one location prohibits
spike generation at other locations. This implies that these coordinates
are highly statistically dependent.

Second, we can compute the true entropy $H(\bX)$ for the spike process
unlike other complicated stochastic processes.
Since the spike process selects one possible vector from the standard basis
of $\R^\ndim$ with uniform probability $1/\ndim$, 
the true entropy $H(\bX)$ is clearly $\log \ndim$.
This is one of the rare cases where we know the true high-dimensional
entropy of the process.

\subsection{The LSDB among the Haar-Walsh Dictionary}
Our first theorem specifies the LSDB selected from the well-known Haar-Walsh 
dictionary, a subset of $\OG(\ndim)$.  This dictionary contains a large
number of orthonormal bases (in fact, more than $2^{\ndim/2}$ bases)
including the standard basis, the Haar basis (consists of dyadically-scaled and
shifted versions of boxcar functions), and the Walsh basis (consisting of 
square waves).  Because the basis vectors in this dictionary are all piecewise
constant (except the standard basis vectors), they are often used to analyze
and compress discontinuous or blocky signals such as acoustic impedance 
profiles of subsurface structure. 
See Wickerhauser (1994), Saito (2000), and Section~\ref{sec:review-Haar-Walsh}
of this paper for the details of this dictionary.
\begin{aismth}
\label{thm:spike-Haar-Walsh}
Suppose we restrict our search of the bases within the Haar-Walsh
dictionary.  Then, the LSDB is:
\begin{itemize}
\item the standard basis if $\ndim > 4$; and
\item the Walsh basis if $\ndim=2$ or $4$.
\end{itemize}
Moreover, the true independence can be achieved only for $\ndim=2$.
Note that $\ndim$ is always a dyadic number in this dictionary.
\end{aismth}

\subsection{The LSDB among $\OG(\ndim)$}
It is curious what happens if we do not restrict ourselves to
the Haar-Walsh dictionary.  Then, we have the following theorem.
\begin{aismth}
\label{thm:spike-OG}
The LSDB among $\OG(\ndim)$ is the following:
\begin{itemize}
\item for $\ndim \geq 5$, either the standard basis or the basis whose matrix
representation is
\begin{equation}
\label{eq:lsdb-OG5}
B_{\OG(\ndim)} = \frac{1}{\ndim} \left[ \begin{array}{ccccc}
 \ndim-2 &   -2   & \cdots &   -2   &  -2 \\
   -2    & \ndim-2& \ddots &        &  -2 \\
\vdots   & \ddots & \ddots & \ddots & \vdots \\
   -2    &        & \ddots & \ndim-2&  -2 \\
   -2    &   -2   & \cdots &   -2   & \ndim-2
\end{array} \right];
\end{equation}
\item for $\ndim = 4$, the Walsh basis, i.e.,
$B_{\OG(4)} = \frac{1}{2}\left[ \begin{array}{cccc}
1 & 1 & 1 & 1   \\
1 & 1 & -1 & -1 \\
1 & -1 & 1 & -1 \\
1 & -1 & -1 & 1
\end{array}\right]$;

\item for $\ndim=3$, $B_{\OG(3)}=\left[ \begin{array}{ccc}
\frac{1}{\sqrt 3} & \frac{1}{\sqrt 6} & \frac{1}{\sqrt 2} \\
\frac{1}{\sqrt 3} & \frac{1}{\sqrt 6} & \frac{-1}{\sqrt 2} \\
\frac{1}{\sqrt 3} & \frac{-2}{\sqrt 6} & 0
\end{array}\right]$; and

\item for $\ndim=2$, $B_{\OG(2)}=\frac{1}{\sqrt{2}}\left[ \begin{array}{cccc}
1 & 1    \\
1 & -1
\end{array}\right]$,
and this is the only case where the true independence is
achieved.
\end{itemize}
\end{aismth}
\begin{aismrm}
There is an important geometric interpretation of \eqref{eq:lsdb-OG5}.
This matrix can also be written as: 
\bdm
I_\ndim - 2 \frac{\bone_\ndim}{\sqrt{\ndim}}\frac{\bone_\ndim^T}{\sqrt{\ndim}}.
\edm
In other words, this matrix represents the \emph{Householder reflection} with
respect to the hyperplane $\{ \by \in \R^\ndim \, | \, \sum_{i=0}^\ndim y_i = 0 \}$ whose unit normal vector is $\bone_\ndim/\sqrt{\ndim}$.
\end{aismrm}

\subsection{The LSDB among $\GL(\ndim,\R)$}
Before discussing the LSDB among a larger class of bases, 
let us remark an important specifics for a discrete stochastic process.

Let $\bX$ be a random vector obeying a discrete stochastic process with
a probability mass function (pmf) $f_\bX$.  This means that there are only
finite number of possible values (or states) $\bX$ can take.  
Clearly the spike process is a discrete process since the only possible values
are $\{\be_1,\ldots,\be_\ndim\}$, the standard basis vectors.
Then, for any invertible transformation $B \in \GL(n,\R)$ with
$\bY=B^{-1}\bX$, be it orthonormal or not, the total entropy 
of the process before and after the transformation is exactly the same.
Indeed, in the definition of discrete Shannon entropy, 
$-\sum_{j} p_j \log p_j$, 
the values that the random variable takes are of no importance; 
only the number of possible values the random variable can take and its pmf
matter. 
In our case, it is clear that the events $ \{\bX=\ba_{i}\}$ and 
$\{\bY=\bb_{i}\}$ where $\bb_i=B^{-1}\ba_i$ are equivalent; 
otherwise the transformation would not be invertible. 
This shows that the corresponding probabilities are equal:
\bdm
\Prob\{\bX=\ba_{i}\}=\Prob\{\bY=\bb_{i}\}.
\edm
Therefore, considering the expression of discrete entropy, this proves that
\bdm
H(\bY)=H(\bX),
\edm
as long as the transformation matrix belongs to $\GL(\ndim,\R)$.
Note that for the continuous case, this is only true if $B \in \SL^\pm(\ndim,\R)$.
Therefore, for a discrete stochastic process like the spike process,
the LSDB among $\GL(\ndim,\R)$ can be selected by just minimizing the sum of
the coordinate-wise entropy as \eqref{eq:min-c-lsdb} as if 
$\Dict=\SL^\pm(\ndim,\R)$. 
In other words, there is no important distinction
in the LSDB selection from $\GL(\ndim,\R)$ and from $\SL^\pm(\ndim,\R)$ for
discrete stochastic processes.
Therefore, we do not have to treat these two cases separately.

\begin{aismth}
\label{thm:spike-GL}
The LSDB among $\GL(\ndim,\R)$ with $\ndim > 2$ is
the following basis pair (for analysis and synthesis respectively):
\begin{equation}
\label{eq:spike-LSDB-SL-analysis}
\begin{hide}
B_{\GL(\ndim,\R)}^{-1} = \left[ \begin{array}{cccccc}
a & a & \cdots & \cdots & \cdots &a \\
b_2 & c_2 & b_2 & \cdots & \cdots & b_2 \\
b_3 & b_3 & c_3 &  b_3  & \cdots & b_3 \\
\vdots & \vdots& \ddots & \ddots  & \ddots & \vdots \\
b_{\ndim-1} & b_{\ndim-1} & \cdots & b_{\ndim-1} & c_{\ndim-1} & b_{\ndim-1}\\
b_\ndim & b_\ndim & \cdots & \cdots & b_\ndim & c_\ndim
\end{array} \right],
\end{hide}
B_{\GL(\ndim,\R)}^{-1} = \left[ \begin{array}{ccccccc}
a & a & \cdots & \cdots & \cdots & \cdots & a \\
b_2 & c_2 & b_2 & \cdots & \cdots & \cdots & b_2 \\
b_3 & b_3 & c_3 &  b_3  & \cdots  & \cdots & b_3 \\
\vdots & \vdots&  & \ddots &      &        & \vdots \\
\vdots & \vdots&  &        & \ddots &     & \vdots \\
b_{\ndim-1} & \cdots & \cdots & \cdots & b_{\ndim-1} & c_{\ndim-1} & b_{\ndim-1}\\
b_\ndim & \cdots & \cdots & \cdots & \cdots & b_\ndim & c_\ndim
\end{array} \right],
\end{equation}
where $a$, $b_k$, $c_k$ are arbitrary real-valued constants satisfying
$a \neq 0$, $b_k \neq c_k$, $k=2,\ldots, \ndim$.
\begin{equation}
\label{eq:spike-LSDB-SL-synthesis}
B_{\GL(\ndim,\R)} = \left[ \begin{array}{ccccc}
\(1 + \sum_{k=2}^\ndim b_kd_k \)/a & -d_2 & -d_3 & \cdots & -d_\ndim \\
-b_2 d_2/a & d_2 & 0 & \cdots & 0 \\
-b_3 d_3/a & 0 & d_3 & \ddots & \vdots \\
\vdots & \vdots & \ddots & \ddots  & 0 \\
-b_\ndim d_\ndim/a & 0 & \cdots & 0 & d_\ndim
\end{array} \right],
\end{equation}
where $d_k = 1/(c_k-b_k)$, $k=2,\ldots,\ndim$.

If we restrict ourselves to $\Dict=\SL^\pm(\ndim,\R)$, then
the parameter $a$ must satisfy:
\bdm
a = \pm \prod_{k=2}^\ndim (c_k - b_k)^{-1}.
\edm 
\end{aismth}

\begin{aismrm}
\label{rem:lsdb-dense}
The LSDB such as \eqref{eq:lsdb-OG5} and the LSDB pair 
\eqref{eq:spike-LSDB-SL-analysis}, \eqref{eq:spike-LSDB-SL-synthesis}
provide us with further insight into the difference
between sparsity and statistical independence.
In the case of \eqref{eq:lsdb-OG5}, this is the LSDB, yet does not 
sparsify the spike process at all. In fact, these coordinates are completely
dense, i.e., $\Cost_0 = \ndim$.  We can also show that the sparsity measure
$\Cost_p$ gets worse as $\ndim \to \infty$.
More precisely, we have the following proposition.
\begin{aismpr}
\label{prop:spike-cp-OG}
\bdm
\lim_{\ndim \to \infty} \Cost_p\(B_{\OG(\ndim)} \, | \, \bX \) = \left\{
\begin{array}{ll}
\infty & \quad \mbox{if $0 \leq p < 1$};\\
3 & \quad \mbox{if $p=1$}.
\end{array}
\right.
\edm
\end{aismpr}
It is interesting to note that this LSDB approaches to the standard
basis as $\ndim \to \infty$.  This also implies that
\bdm
\lim_{\ndim \to \infty} \Cost_p\(B_{\OG(\ndim)} \, | \, \bX \) \neq \Cost_p\(\lim_{\ndim \to \infty} B_{\OG(\ndim)} \, | \, \bX \).
\edm

As for the analysis LSDB \eqref{eq:spike-LSDB-SL-analysis}, 
the ability to sparsify the spike process depends on the values of $b_k$ 
and $c_k$.
Since the parameters $a$, $b_k$ and $c_k$ are arbitrary as long as 
$a \neq 0$ and $b_k \neq c_k$,
let us put $a=1$, $b_k=0$, $c_k=1$, for $k=2,\ldots,\ndim$.
Then we get the following specific LSDB pair:
\bdm
\label{eq:spike-LSDB-SL-01}
B^{-1}_{\GL(\ndim,\R)} = \left[ \begin{array}{cccc}
1 & 1 & \cdots & 1 \\
0 & & &  \\
\vdots & & I_{\ndim-1} & \\
0 & & & \end{array} \right],
\quad
B_{\GL(\ndim,\R)} = \left[ \begin{array}{cccc}
1 & -1 & \cdots & -1 \\
0 & & & \\
\vdots & & I_{\ndim-1} & \\
0 & & & \end{array} \right].
\edm
This analysis LSDB provides us with a sparse representation for the spike
process (though this is clearly not better than the standard basis). 
For $\bY=B_{\GL(\ndim,\R)}^{-1}\bX$,
\bdm
\Cost_0 = E \left[ \| \bY \|_0 \right]= \frac{1}{\ndim} \times 1 +\frac{\ndim-1}{\ndim} \times 2=2-\frac{1}{\ndim}.
\edm
Now, let us take $a=1$, $b_k=1$, $c_k=2$ for $k=2,\ldots,\ndim$ in \eqref{eq:spike-LSDB-SL-analysis} and \eqref{eq:spike-LSDB-SL-synthesis}. Then we get
\begin{equation}
\label{eq:spike-LSDB-SL-12}
\begin{hide} 
B^{-1}_{\GL(\ndim,\R)} = \left[ \begin{array}{ccccc}
1 & 1 & 1 &  \cdots & 1 \\
1 & 2 & 1 &  \cdots & 1 \\
1 & 1 & 2 &         & 1 \\
\vdots & \vdots& & \ddots  & \vdots \\
1 & 1 & 1 &         & 2 
\end{array} \right],
\quad
B_{\GL(\ndim,\R)} = \left[ \begin{array}{ccccc}
\ndim & -1 & -1 & \cdots & -1 \\
-1    &  1 & 0  & \cdots &  0 \\
-1    &  0 & 1  &        &  0 \\
\vdots & \vdots& & \ddots  & \vdots \\
-1    &  0 & 0   &        & 1
\end{array} \right].
\end{hide}
\begin{hide} 
B^{-1}_{\GL(\ndim,\R)} = \left[ \begin{array}{ccccc}
1 & 1 & \cdots & \cdots & 1 \\
1 & 2 & \ddots &        & \vdots \\
\vdots & \ddots & \ddots & \ddots & \vdots \\
\vdots &        & \ddots  & 2 & 1 \\
1 & \cdots & \cdots & 1   & 2 
\end{array} \right],
\quad
B_{\GL(\ndim,\R)} = \left[ \begin{array}{ccccc}
\ndim & -1 & \cdots & \cdots & -1 \\
-1    &  1 & 0  & \cdots &  0 \\
\vdots & 0 & \ddots  & \ddots  &  \vdots \\
\vdots & \vdots& \ddots & \ddots  & 0 \\
-1    &  0 & \cdots   &     0   & 1
\end{array} \right].
\end{hide}
B^{-1}_{\GL(\ndim,\R)} = \left[ \begin{array}{cccc}
1 & 1 & \cdots & 1 \\
1 & 2 & \ddots & \vdots \\
\vdots & \ddots & \ddots & 1 \\
1 & \cdots & 1 & 2 
\end{array} \right],
\quad
B_{\GL(\ndim,\R)} = \left[ \begin{array}{cccc}
\ndim  & -1 & \cdots &  -1 \\
-1     &   &   &   \\
\vdots &   & I_{\ndim-1}  &   \\
-1     &   &   &   \\ 
\end{array} \right].
\end{equation}
The spike process under this analysis basis is completely dense,
i.e., $\Cost_0 = \ndim$. Yet this is still the LSDB.
\end{aismrm}

Finally, from Theorems~\ref{thm:spike-OG} and \ref{thm:spike-GL}, we can prove
the following corollary:
\begin{aismco}
\label{cor:spike-GL}
There is no invertible linear transformation providing the 
statistically independent coordinates for the spike process for $\ndim > 2$.
In fact, the mutual information $I\(B_{\OG(\ndim)}^T \bX\)$ and
$I\(B_{\GL(\ndim,\R)}^{-1} \bX\)$ are monotonically increasing as a function
of $\ndim$, and both approaches to $\log \e \approx 1.4427$ as $\ndim \to \infty$.
\end{aismco}
\begin{aismrm}
\label{rem:spike-convert}
Although the spike process is very simple, we have the following 
interpretation. 
Consider a stochastic process generating a basis vector randomly at a time 
selected from some orthonormal basis.  Then, both that basis itself is the BSB
and the LSDB among $\OG(\ndim)$.  Theorem~\ref{thm:spike-OG} claims that 
once we transform the data to the spikes, one cannot do any better than that 
both in sparsity and independence within $\OG(\ndim)$.
Of course, if one extends the search to nonlinear transformations,
then it becomes a different story.  We refer the reader to our
recent articles Lin et al.\ (2000, 2001) for the details of a nonlinear 
algorithm.
\end{aismrm}

\section{Proofs of Propositions and Theorems}
\label{sec:proofs}
\subsection{Proof of Proposition~\ref{prop:spike-KLB}}
{\sc Proof.}
Let $\bX=(X_{1},X_{2},\ldots,X_{\ndim})^T$ be a random vector generated by this
process. 
For each of its realizations,  a randomly chosen coordinate among these 
$\ndim$ positions takes the value $1$, while the others take the value $0$.
Hence each $X_i$, $i=1,\ldots,\ndim$, takes the values $1$ with probability 
$1/\ndim$ and the value $0$ with probability $1-1/\ndim$.
Let us calculate the covariance of these variables.  First, we have:
\begin{eqnarray*}
E(X_{i})&=&\frac{1}{\ndim} \times 1+ \( 1-\frac{1}{\ndim} \) \times 0=\frac{1}{\ndim}  \quad \mbox{for $i=1, \ldots, \ndim$} \\
E(X_{i}X_{j}) &=& 
  \left\{ \begin{array}{ll}
	E(X_{i}^2)=E(X_{i})\;& \mbox{if $i=j$}; \\
	0 \;& \mbox{if $i \neq j$,}
          \end{array}
  \right.
\end{eqnarray*}
since one of these two variables will always take the value $0$.
Let $R=(R_{ij})$ be the covariance matrix of this process.
Then, we have:
\bdm
R_{ij}=E(X_{i}X_{j})-E(X_{i})E(X_{j})= \frac{1}{\ndim} \delta_{ij} - \frac{1}{\ndim^2}
\edm
We know that a basis is a Karhunen-Lo\`eve basis if and only if it is 
orthonormal and diagonalizes the covariance matrix. 
Thus, we will now calculate the eigenvalue decomposition of the covariance 
matrix $R=\frac{1}{\ndim}I_{\ndim}-\frac{1}{\ndim^2}J_{\ndim}$,
where $I_{\ndim}$ is the identity matrix of size $\ndim \times \ndim$, 
and $J_{\ndim}$ is an $\ndim \times \ndim$ matrix with each entry taking
the value $1$. 

We now need to calculate the determinant:
\begin{displaymath} 
 P_{R}(\lambda) \define \det(\lambda I_\ndim - R) = 
   \left|
	 \begin{array}{cccc}
	\lambda-\frac{1}{\ndim}+\frac{1}{\ndim^2} & \frac{1}{\ndim^2} &\ldots & \frac{1}{\ndim^2} \\
	\frac{1}{\ndim^2} & \ddots & \ddots &  \vdots \\
	\vdots &  \ddots &\ddots & \frac{1}{\ndim^2}         \\
	\frac{1}{\ndim^2} & \ldots & \frac{1}{\ndim^2} & \lambda-\frac{1}{\ndim}+\frac{1}{\ndim^2}
	\end{array}
   \right|,
\end{displaymath}
which is of the generic form:
\begin{displaymath} 
 \Delta(a,b) \define \left|
	\begin{array}{cccc}
   	a+b & b &\ldots & b\\
	b & a+b & \ddots & \vdots \\
	\vdots &\ddots &\ddots & b \\
	b & \ldots & b & a+b 
	\end{array}
   \right|,
\end{displaymath}
with the values $a=\lambda-1/\ndim$  and $b=1/\ndim^2$.
This is calculated by subtracting the last row from all the others, 
and then adding all $\ndim-1$ columns to the last one.

\begin{equation}
\label{eq:delta-a-b}
 \Delta(a,b)=\left|
	\begin{array}{ccccc}
   	a & 0 &\ldots &0 & -a\\
	0 & a & \ddots & \vdots &\vdots \\
	\vdots &\ddots &\ddots&0& \vdots \\
	0&\ldots&0&a &-a \\
	b & \ldots & \ldots & b & a+b 
	\end{array}
   \right| =\left|
	\begin{array}{ccccc}
   	a & 0 &\ldots &0 & 0\\
	0 & a & \ddots & \vdots &\vdots \\
	\vdots &\ddots &\ddots&0& \vdots \\
	0 &\ldots&0&a &0 \\
	b & \ldots & \ldots & b & a+\ndim b 
	\end{array}
   \right| =a^{\ndim-1}(a+\ndim b).
\end{equation}
Putting $a=\lambda-1/\ndim$ and $b=1/\ndim^2$,
we have the characteristic polynomial $P_R$ of $R$ as
$P_{R}(\lambda)= \lambda (\lambda-1/\ndim)^{\ndim-1}$.
Hence, the eigenvalues of $R$ are $\lambda=0$ or $1/\ndim$.


It is now obvious that the vector $\bone_\ndim=(1,\ldots,1)^T$
is an eigenvector for $R$ associated with the eigenvalue $0$, i.e.,
$\bone_\ndim \in \ker R$.
Indeed, we have
\bdm
R\bone_\ndim=\(\frac{1}{\ndim}I_{\ndim}-\frac{1}{\ndim^2}J_{\ndim}\)\bone_\ndim=\frac{1}{\ndim}\ \bone_\ndim-\frac{1}{\ndim^2}\ \ndim \bone_\ndim=0.
\edm
Since $\dim \ker R =1$, $\ker R$ is the one-dimensional subspace spanned by
$\bone_\ndim$ . Considering that $R$ is symmetric and only has two distinct 
eigenvalues, we know that the eigenspace associated to the eigenvalue 
$1/\ndim$ is orthogonal to $\ker R$, which is the hyperplane
$\{ \by \in \R^\ndim \, | \, \sum_{i=1}^{\ndim}y_{i}=0 \}$.
Therefore, the orthogonal bases that diagonalize $R$ are the bases formed by 
the adjunction of $\bone_\ndim$ to any orthogonal basis of $\ker R^\perp$.
The Walsh basis, which consists of oscillating square waves, 
is such a basis, although it is just one among many.
{\hfill $\Box$}

\subsection{Proof of Proposition~\ref{prop:spike-sparse}}
{\sc Proof.}
The case $\Dict=\OG(\ndim)$ is obvious as discussed before this proposition.
Therefore, we first prove the case $\Dict=\GL(\ndim,\R)$.
To maximize the sparsity, it is clear that the transformation matrix
must be diagonal (modulo permutations and sign flips), i.e., 
$B_p = \diag(a_1,\ldots,a_\ndim)$ with $a_k \neq 0$, $k=1,\ldots,\ndim$.
The sparsity cost $\Cost_p$ defined in \eqref{eq:cp} can be computed
and bounded in this case as follows:
\bdm
\Cost_p (B_p \, | \, \bX) = E \| \bY \|^p_p = \frac{1}{\ndim} \sum_{k=1}^\ndim | a_k |^p
\geq | a |^p,
\edm
where $| a | = \min \left\{ | a_1 |, \ldots, | a_\ndim | \right\}$.
This lower bound is achieved when $B_p= a I_\ndim$, i.e., a nonzero
constant times the standard basis.
Now, if $\Dict=\SL^\pm(\ndim,\R)$, then this constant $a$ must
be either $1$ or $-1$ since $\det(B_p) = a^\ndim=\pm 1$ and $a \in \R$.
{\hfill $\Box$}

\subsection{A Brief Review of the Haar-Walsh Dictionary and the Best-Basis Algorithm}
\label{sec:review-Haar-Walsh}
Before proceeding to the proof of Theorem \ref{thm:spike-Haar-Walsh},
let us first review the Haar-Walsh dictionary and define some necessary
quantities. 

Let $\ndim$  be a positive dyadic integer, i.e., $\ndim=2^\ndimlev$ for some 
$\ndimlev \in \N$.
An input vector $\bx=(x_1,\ldots,x_\ndim)^T$, viewed as a digital signal
sampled on a regular grid in time, is first decomposed into low and high 
frequency bands by the convolution-subsampling operations on the discrete
time domain with the pair consisting of a ``lowpass'' filter 
$\{h_\ell\}_{\ell=1}^{\LF}$ and a ``highpass'' filter $\{g_\ell\}_{\ell=1}^{\LF}$.
Let $H$ and $G$ be the convolution-subsampling operators using these filters
which are defined as:
\bdm
(H\bx)_k = \sum_{\ell=1}^{\LF} h_\ell x_{\ell+2(k-1)}, \quad
(G\bx)_k = \sum_{\ell=1}^{\LF} g_\ell x_{\ell+2(k-1)},
\quad k=1,\ldots,\ndim.  
\edm
We assume the periodic boundary condition 
on $\bx$ (whose period is $\ndim$).  Hence, the filtered sequences 
$H\bx$ and $G\bx$ are also periodic with period $\ndim/2$.
Their adjoint operations (i.e., upsampling-anticonvolution) $H^*$ and $G^*$
are defined as 
\bdm
(H^*\bx)_k = \sum_{1 \leq k-2(\ell-1) \leq \LF} h_{k-2(\ell-1)} x_\ell, \quad
(G^*\bx)_k = \sum_{1 \leq k-2(\ell-1) \leq \LF} g_{k-2(\ell-1)} x_\ell,
\quad k=1,\ldots,2\ndim.
\edm
The filter $H$ and $G$ are called conjugate mirror filters (CMF's) if 
they satisfy the following orthogonality (or perfect reconstruction) 
conditions: 
\bdm
HG^* = GH^* = 0 \quad \mbox{and} \quad H^* H + G^* G = I,
\edm
where $I$ is the identity operator. 
Various design criteria (concerning regularity, symmetry etc.)\ on the lowpass
filter coefficients $\{h_\ell\}$ can be found in Daubechies (1992).
The Haar-Walsh dictionary uses the filter pair with the shortest length 
($\LF=2$) and $h_1 = h_2 = 1/\sqrt{2}$.
Once $\{h_\ell\}$ is fixed, the filter $G$ is obtained by setting 
$g_\ell = (-1)^{\ell-1} h_{\LF-\ell+1}$.
This decomposition process is iterated on both the low and high frequency
components.
The first level decomposition generates two subsequences $H\bx$ and $G\bx$
each of which has length $\ndim/2$.
In the case of the Haar-Walsh dictionary, these subsequences are:
\bdm
H\bx =\( \frac{x_{1}+x_{2}}{\sqrt2},  \ldots, \frac{x_{\ndim-1}+x_{\ndim}}{\sqrt2} \)^T \quad \mbox{and}  \quad G\bx = \( \frac{x_{1}-x_{2}}{\sqrt2},  \ldots, \frac{x_{\ndim-1}-x_{\ndim}}{\sqrt2} \)^T.
\edm
The second level generates four subsequences,
$H^2\bx$, $GH\bx$, $HG\bx$, $G^2\bx$, each of which is of length $\ndim/4$.
If we repeat this process for $k$ times ($k=0,1,\ldots,\maxlev \leq \ndimlev$),
then at the $k$th level,
$2^k$ subsequences $H^k\bx$, $GH^{k-1}\bx$, $\ldots$, $G^{k-1}H\bx$, $G^k\bx$,
each of which is of length $2^{\ndimlev-k}$, are generated.
As a whole there are $(k+1) \ndim$ expansion coefficients
(including the original components of $\bx$).
One can iterate this procedure and stop at the $\maxlev$th level, where
$\maxlev \leq \ndimlev$.
These coefficients are naturally organized in the binary tree structure
as shown in Figure~\ref{fig:wptable}.
\begin{figure}
\begin{center}
\begin{hide}
\begin{tabular}{l *{16}{|c| p{1.5mm}}} \cline{2-17}
$k=0$ & \multicolumn{16}{|c|}{+} \\  \cline{2-17}
$k=1$ & \multicolumn{8}{|c|}{+} & \multicolumn{8}{|c|}{-} \\ \cline{2-17}
$k=2$ & \multicolumn{4}{|c|}{+} & \multicolumn{4}{|c|}{-} &
\multicolumn{4}{|c|}{-} & \multicolumn{4}{|c|}{-} \\ \cline{2-17}
\vdots  & \multicolumn{16}{|c|}{\ \ \ \ \ldots \ \ \ \ } \\ \cline{2-17}
$k=\maxlev$ & + & - & - & - & - & - & - & - & - & - & - & - & - & - & - & -  \\ \cline{2-17}
\end{tabular} 
\end{hide}
\begin{tabular}{c *{16}{|c} |} \cline{2-17}
$k=0$ & \multicolumn{16}{|c|}{+} \\  \cline{2-17}
$k=1$ & \multicolumn{8}{|c|}{+} & \multicolumn{8}{|c|}{-} \\ \cline{2-17}
$k=2$ & \multicolumn{4}{|c|}{+} & \multicolumn{4}{|c|}{-} &
\multicolumn{4}{|c|}{-} & \multicolumn{4}{|c|}{-} \\ \cline{2-17}
\vdots  & \multicolumn{16}{|c|}{ $\cdots$ } \\ \cline{2-17}
$k=\maxlev$ & + & - & - & - & - & - & - & - & - & - & - & - & - & - & - & - \\ \cline{2-17} 
\end{tabular} 
\vspace{1em}
\caption{A table of dictionary coefficients are organized as the binary
tree structured table.}
\label{fig:wptable}
\end{center}
\end{figure}
For future reference, we refer to the tree with $\maxlev=\ndimlev$
as the \emph{maximal-depth} tree or the \emph{full} tree.
Because of the perfect reconstruction condition on $H$ and $G$,
each decomposition step is also interpreted as a decomposition of the vector
space into mutually orthogonal subspaces.  Let $\subsp_{0,0}$ denote
the $\ndim$-dimensional Euclidean space $\R^\ndim$ spanned by the standard 
basis vectors.  Hence, an input vector of length $\ndim$ is an element
of $\subsp_{0,0}$.
Let $\subsp_{1,0}$ and $\subsp_{1,1}$ be mutually
orthogonal subspaces generated by the application of the operators $H$ and 
$G$ respectively to the parent space $\subsp_{0,0}$.
Then, in general, the $k$th step of the decomposition process 
($k=0,\ldots,\maxlev$) can be written as
\bdm
\subsp_{k,\ell} = \subsp_{k+1,2\ell} \oplus \subsp_{k+1,2\ell+1} \quad 
\ell=0,\ldots,2^k-1.
\edm
It is clear that $\dim \subsp_{k,\cdot} = 2^{\ndimlev-k}$. 
For each subspace $\subsp_{k,\ell}$, we associate the basis vectors
$\bw_{k,\ell,m} \in \R^\ndim$, $m=0,\ldots,2^{\ndimlev-k}-1$ which span
this subspace.
The vector $\bw_{k,\ell,m}$ is roughly centered at $2^km$,
has length of support $\approx 2^k$, and oscillates $\approx \ell$ times.
Note that for $k=0$, we have the standard basis of $\R^\ndim$.
The expansion coefficients computed by the convolution-subsampling operations
can be viewed as the inner products between the input vector and these basis
vectors although we never need to compute these inner products explicitly.
Clearly, we have a redundant set of subspaces in the binary tree.
In fact, it is easily proved that there are more than
$2^{2^{\maxlev-1}}$ possible orthonormal bases in this binary tree; see e.g.
Wickerhauser (1994) for the details.  
Because of this abundance of the bases, such a binary tree of subspaces
(or basis vectors) is called a wavelet packet dictionary for general CMF's and 
\emph{the Haar-Walsh dictionary} if $\LF=2$ and $h_1=h_2=1/\sqrt{2}$.
Now an important question is how to select the best coordinate system
efficiently for the problem at hand from this dictionary.

The ``best-basis'' algorithm of Coifman and Wickerhauser (1992)
first expands an input vector into a specified basis dictionary.
Then a complete basis called a \emph{best basis} (BB) which minimizes a certain
cost function (such as the sparsity cost $\Cost_p$ \eqref{eq:cp} or 
the statistical dependence cost $\Cost_H$ \eqref{eq:c-lsdb}; see also 
Saito (2000) for a variety of cost functions for different problems 
such as classification and regression) is searched in this binary tree using
the divide-and-conquer algorithm.
More precisely,
let $B_{k,\ell}$ denote a set of basis vectors belonging to the subspace
$\subsp_{k,\ell}$ arranged as a matrix
\begin{equation}
\label{eq:basis-matrix-jk}
B_{k,\ell} = (\bw_{k,\ell,0},\ldots,\bw_{k,\ell,2^{\ndimlev-k}-1}).
\end{equation}
Now let $A_{k,\ell}$ be the best basis for the input signal $\bx$ restricted
to the span of $B_{k,\ell}$ and let $\Cost$ be a cost function measuring 
the deficiency of the nodes (subspaces) such as $\Cost_p$ or $\Cost_H$.
The following best-basis algorithm ``prunes'' this binary tree
by comparing the cost of each parent node and its two children nodes:

\noindent
Given an input vector $\bx \in \R^\ndim$,\\[-20pt]
\begin{description}
\item[Step 0:] Choose a basis dictionary $\Dict$,
specify the maximum depth of decomposition $\maxlev$, and an information cost $\Cost$.
\item[Step 1:] Expand $\bx$ into the dictionary $\Dict$
and obtain coefficients $\{B_{k,\ell}^T \bx\}_{0 \leq k \leq \maxlev,\ 0 \leq \ell  \leq 2^k-1}$.
\item[Step 2:] Set $A_{\maxlev,\ell} = B_{\maxlev,\ell}$ for $\ell=0,\ldots,2^\maxlev-1$.
\item[Step 3:] Determine the best subspace $A_{k,\ell}$ in the bottom-up
manner, i.e., for $k=\maxlev-1,\ldots,0, \ \ell = 0,\ldots,2^k-1$, by
\begin{equation}
\label{eq:bb-recurse}
A_{k,\ell} = \left\{ \begin{array}{ll}
	B_{k,\ell} & \mbox{if \ $\Cost(B_{k,\ell}^T \bx) \leq \Cost(A_{k+1,2\ell}^T\bx \cup  A_{k+1,2\ell+1}^T\bx)$},\\
	A_{k+1,2\ell} \oplus A_{k+1,2\ell+1} & \mbox{otherwise.}
\end{array}
\right.
\end{equation}
\end{description}
This algorithm becomes fast if the cost function $\Cost$ is \emph{additive},
i.e., 
$\Cost(\bzero) = 0$ and $\Cost(\bx) = \sum_{i} \Cost(x_i)$.
Both $\Cost_p$ of \eqref{eq:cp} and $\Cost_H$ of \eqref{eq:c-lsdb} are
clearly additive.
If $\Cost$ is additive, then in \eqref{eq:bb-recurse} we have
\bdm
\Cost(A_{k+1,2\ell}^T \bx \cup A_{k+1,2\ell+1}^T \bx) = \Cost(A_{k+1,2\ell}^T \bx) + \Cost(A_{k+1,2\ell+1}^T \bx).
\edm
This implies that a simple addition suffices instead of computing the cost of
union of the nodes.

Coming back to the Haar-Walsh case, we need a few more definitions for
the proof of Theorem~\ref{thm:spike-Haar-Walsh}.
At each level of the decomposition, the leftmost node (or box) representing 
the coefficients $H^k \bx$ is marked by $+$ in Figure~\ref{fig:wptable}.
This node also corresponds to the subspace $\subsp_{k,0}$.
Clearly, each coefficient in this node must be of the form
\begin{equation}
\label{eq:pos-coef}
\frac{1}{\sqrt{2^{k}}}\ \(  x_{\sigma (1)}+\cdots+x_{\sigma (2^k)} \)
\end{equation}
where $\sigma$ is a permutation of $\{1, \ldots, \ndim\}$.
We call these nodes and the corresponding
coefficients the \emph{positive node} and the \emph{positive coefficients},
respectively.
All the other nodes marked by $-$ sign at the $k$th level corresponding to
the subspaces $\subsp_{k,\ell}$, $\ell \neq 0$, contain
coefficients of the form:
\begin{equation}
\label{eq:neg-coef}
\frac{1} {\sqrt{2^{k}}}\ \(
 x_{\sigma (1)}+\ldots+x_{\sigma (2^{k-1})}-x_{\sigma (2^{k-1}+1)}-\ldots-x_{\sigma (2^k)}  \).
\end{equation}
These nodes and coefficients are referred to as \emph{negative nodes} and 
\emph{negative coefficients}, respectively.
We note that any descendant node of a negative node must be negative.
In fact, only the left child node of a positive node can be positive.

\subsection{Proof of Theorem~\ref{thm:spike-Haar-Walsh}}
\begin{figure}
\epsfxsize=0.6\textwidth
\centerline{\epsfbox{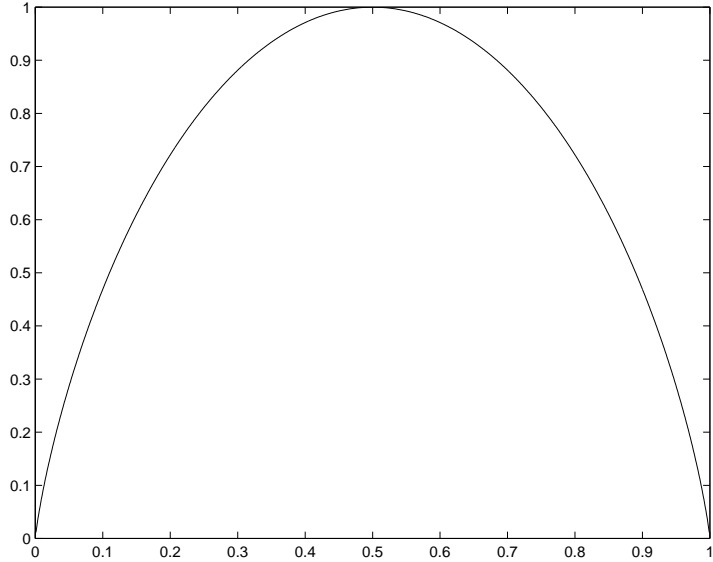}}
\caption{A plot of $f:x \rightarrow - \left[ x \log{x} + (1-x)\log(1-x) \right] $.}
\label{fig:f}
\end{figure}

{\sc Proof.}
Let us consider the positive coefficients. The $k$th-level positive node
contains $\ndim/2^k$ coefficients each of which is generated by \eqref{eq:pos-coef}, which in the case of the spike process can take only the following
values:
\begin{itemize}
\item $+1/\sqrt{2^k}$ with probability $2^k/\ndim$;
\item $0$ with probability $1-2^k/\ndim$. 
\end{itemize}
Thus the entropy of each coordinate in the $k$th-level positive node can
be computed as
\bdm
h_{+}(k) \define -\( \frac{2^k}{\ndim}\ \log \( \frac{2^k}{\ndim}\)+\(1- \frac{2^k}{\ndim} \) \log\(1-\frac{2^k}{\ndim}\) \)=f\( \frac{2^k}{\ndim} \),
\edm
where
\begin{equation}
\label{eq:f}
f(x) \define -[x\log(x)+(1-x) \log(1-x)],
\end{equation}
which is displayed in Figure~\ref{fig:f}.
The following properties of this function $f$ are basic and will be used 
repeatedly in this paper:
\begin{itemize}
\item For all $x \in [0,1]$, $f(x)\geq 0 $ and $f(x)=0$ if and only if $x=0$
or $x=1$;
\item For all $x \in [ 0,1]$, $f(x)=f(1-x)$;
\item $f$ is increasing on $[0,1/2]$, and decreasing on $[1/2,1]$;
\item $f$  is concave on $[0,1]$.
\end{itemize}

\begin{figure}
\epsfxsize=0.6\textwidth
\centerline{\epsfbox{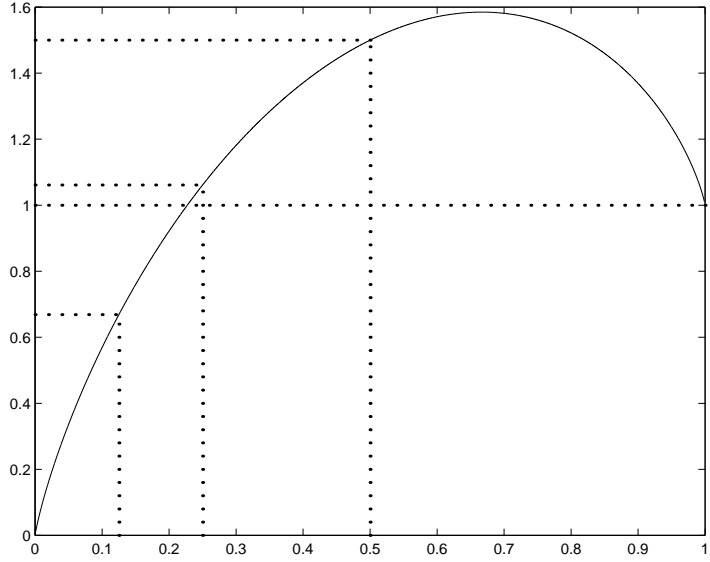}}
\caption{A plot of $g:x \rightarrow - \left[ x \log{\frac{x}{2}} + (1-x)\log(1-x) \right] $.}
\label{fig:g}
\end{figure}

On the other hand, the remaining $\ndim-(\ndim/2^k)$ negative 
coefficients at level $k$ are computed by \eqref{eq:neg-coef},
which can take three different values:
\begin{itemize}
\item $+1/\sqrt{2^k}$ with probability $2^{k-1}/\ndim$;
\item $-1/\sqrt{2^k}$ with probability $2^{k-1}/\ndim$;
\item $0$ with probability $1-2^k/\ndim$. 
\end{itemize} 
Thus the entropy of each negative coordinate of level $k$ is
\bdm
h_{-}(k)\define-\( \frac{2^k}{\ndim} \log\( \frac{2^{k-1}}{\ndim}\)+\(1- \frac{2^k}{\ndim}\)\log\(1-\frac{2^k}{\ndim}\) \)=g\( \frac{2^k}{\ndim} \),
\edm
where
\begin{equation}
\label{eq:g}
g(x) \define -[x \log (x/2)+(1-x) \log (1-x)]=f(x)+x,
\end{equation}
which is plotted in Figure~\ref{fig:g}.

\begin{figure}
\epsfxsize=0.6\textwidth
\centerline{\epsfbox{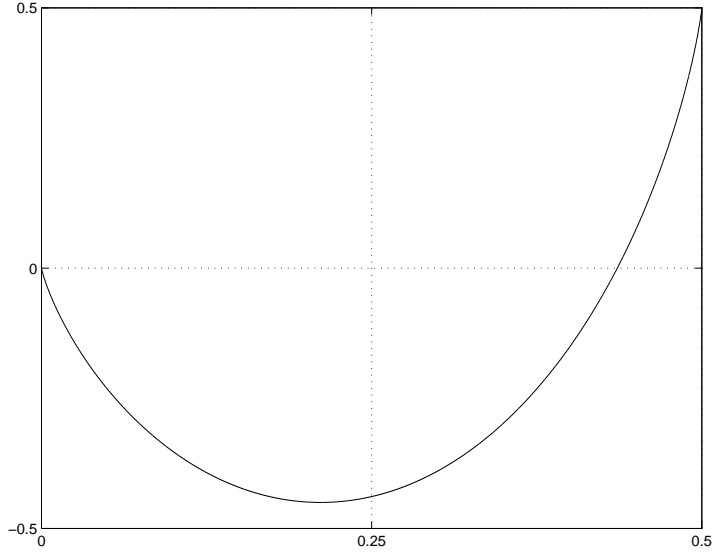}}
\caption{A plot of the function $g(x)-g(2x) $.}
\label{fig:gg2}
\end{figure}

The following lemma is used to compare the entropy cost between
a parent node and its children nodes of the Haar-Walsh dictionary.
\begin{aismle}
\label{lem:h+-}
\begin{eqnarray}
h_{-}(k)& \leq & h_{-}(k+1) \label{eq:neg}   \\
h_{+}(k)& \leq & \frac{1}{2} \left[h_{+}(k+1)+h_{-}(k+1)\right] \label{eq:pos} ,
\end{eqnarray} 
for $k=1, \ldots, \ndimlev-2$.
\end{aismle}
{\sc Proof.}
Using the function $g$ defined in \eqref{eq:g}, we have
$h_{-}(k)- h_{-}(k+1)=g(2^k/\ndim)-g(2^{k+1}/\ndim)$. 
As shown in Figure~\ref{fig:gg2}, the function $g(x)-g(2x)$ is always
negative as long as 
$x=2^k/\ndim \leq 0.43595\cdots$.  Since $\ndim=2^\ndimlev$,
this implies that $k-\ndimlev \leq \log(0.43595) \approx -1.1977$, i.e.,
$k \leq \ndimlev - 2$.
Hence we have proved \eqref{eq:neg}.
To prove \eqref{eq:pos}, we have 
$h_+(k)-\frac{1}{2} [ h_+(k+1)+h_-(k+1) ]=f(2^k/\ndim)-\frac{1}{2} [ f(2^{k+1}/\ndim)+ g(2^{k+1}/\ndim) ]$. However, 
\begin{eqnarray*}
f(x)-\frac{1}{2}  \left[ f(2x)+g(2x) \right]
&=& f(x)-\frac{1}{2} \left[g(2x)-2x+g(2x)\right] \\
&=&(f(x)+x)-g(2x) \\
&=& g(x)-g(2x) < 0,
\end{eqnarray*}
if $x=2^k/\ndim \leq 0.43595$, i.e., $k \leq \ndimlev-2$ as before.
{\hfill $\Box$}\vspace{1em}

Inequality~\eqref{eq:neg} implies that the entropy corresponding to a negative
coordinate at one level is smaller than that of the level below. 
Therefore, a negative parent node has smaller entropy than its two negative
children nodes provided that the children nodes are \emph{not 
the bottom leaves}, i.e., if the maximal decomposition level $\maxlev$
satisfies $\maxlev < \ndimlev$. 
In fact, we have
$h_{-}(\ndimlev-2) \leq  h_{-}(\ndimlev-1)$, but 
$h_{-}(\ndimlev-1) \geq  h_{-}(\ndimlev)$.
This means that starting from the non-maximal depth negative nodes,
the best-basis algorithm always chooses the furthest possible ancestor 
negative nodes.

As for the positive nodes, from \eqref{eq:pos}, we can compare the
total entropy of the positive node at level $k$ with that of the two
children nodes (positive and negative) as follows:
\bdm
\frac{\ndim}{2^k} h_{+}(k) \leq \frac{\ndim}{2^{k+1}} \left[ h_{+}(k+1) + h_{-}(k+1) \right],
\edm
since the parent node contains
$\ndim/2^k$ coordinates and each of the children node has
$\ndim/2^{k+1}$ coordinates.
Therefore, again the parent positive node has smaller entropy than its
two children nodes as long as the tree is of non-maximal depth.

These two facts prove that the best-basis algorithm seeking the minimum
entropy selects the root node, i.e., the LSDB is the standard basis,
if $\maxlev \leq \ndimlev-1$.

Now, we need to consider the case of the maximal-depth tree.
Notice that although \eqref{eq:neg} does not hold for $k=\ndimlev-1$,
the following holds: 
\bdm
h_{-}(\ndimlev-3) \leq  h_{-}(\ndimlev),
\edm
since $g(1) \geq g(1/8)$ (see also Figure~\ref{fig:g}).
This allows us, using the best-basis algorithm, to move up from a pair of
bottom leaves not to their immediate parent but to their ``great-grandfather'',
with decreasing entropy, as long as this great-grandfather is still 
a negative node and $\ndimlev \geq 3$.
We still need to consider what happens if this assumption is false, 
that is, if we have maximal-depth leaves with positive great-grandfather.
The self-similar structure of this tree proves that this problem is 
equivalent to the general problem with $\ndimlev=3$, 
which we shall now discuss. 

\begin{description}
\item[\underline{$\ndimlev=3$ (i.e., $\ndim=8$):}]
\begin{hide}
\begin{center}
\begin{tabular}{l *{8}{|c  }|} \cline{2-9}
k=0 & \multicolumn{8}{|c|}{+} \\  \cline{2-9}
k=1 & \multicolumn{4}{|c|}{+} & \multicolumn{4}{|c|}{-} \\ \cline{2-9}
k=2 & \multicolumn{2}{|c|}{+} & \multicolumn{2}{|c|}{-} &
\multicolumn{2}{|c|}{-} & \multicolumn{2}{|c|}{-} \\ \cline{2-9}
k=3 & + & - & - & - & - & - & - & -  \\ \cline{2-9}
\end{tabular}
\end{center}
\end{hide}
Let us show that whatever the set of coordinates chosen among these, 
the entropy they generate is larger than that of the root node, which is
also positive.  The entropy of the root node is:
$8 \times f(1/8)\simeq 4.34$ bits.

The choice of a basis in this dictionary is equivalent to the choice of 
a binary tree of depth $\maxlev \leq 3$.
This reduces to:
\begin{itemize}
\item the choice of the level of the positive node in the basis, 
which also amounts to the choice of the depth of the leftmost leave of 
the tree.
\item the choice of an orthonormal basis of the subspace orthogonal to 
the chosen positive node.
\end{itemize}
We note that all the negative coordinates of the tree have larger entropy
than those of the bottom leaves: this is derived from 
$g(1) \leq g(1/4) \leq g(1/2)$ (see Figure~\ref{fig:g}). 
Thus, the entropy of any basis with its positive node on level $k$ 
is larger than $2^{\ndimlev-k} \times f(2^{k-\ndimlev})+(n-2^{\ndimlev-k}) \times g(1)$.
 
Then there are three different cases corresponding to the level of 
the positive node:
\begin{itemize}
\item if the positive node is on the bottom level, then we only have one
positive coordinate, and seven negative ones; therefore, the entropy of
any such basis is larger than $f(1)+7 \times g(1)=7$ bits; 
\item if the positive node is on level $k=2$, 
we have two positive coordinates and six negative ones; thus the entropy 
of any such basis is larger than $2 \times f(1/2)+6 \times g(1) \simeq 8$ bits;
\item finally, if the positive node is on level $k=1$, then the entropy of
any such basis is larger than  $4 \times f(1/4) +4 \times g(1) \simeq 7.24$ bits.
\end{itemize}
All these values are larger than the entropy of the root node of this tree, 
namely the standard basis. 
Therefore, the standard basis is the LSDB among the Haar-Walsh dictionary 
for $\ndimlev=3$.

\item[\underline{$\ndimlev \geq 3$ (i.e., $\ndim \geq 8$):}]
What we saw shows that this is also true for any integer $\ndimlev \geq 3$
thanks to the self-similar structure of the binary tree dictionary.
This ends the proof of the first part of the theorem: the standard basis is
the LSDB among the Haar-Walsh dictionary for $\ndimlev \geq 3$, i.e.,
$\ndim \geq 8$.
\end{description}
Therefore, we are left to consider the two special cases $\ndim=2$ and
$\ndim=4$.

\begin{description}
\item[\underline{$\ndim=2$:}] 
In this case, the components of the spike process in
the Walsh basis are truly independent. 
Indeed, the representation of $\bx=(x_1,x_2)^T$ in the Walsh basis is:
$(\frac{x_1+x_2}{\sqrt{2}},\frac{x_1-x_2}{\sqrt{2}})^T$.
The sum of the coordinate-wise entropy of the spike process 
relative to the Walsh basis is $h_+(1)+h_-(1)=f(1)+g(1)=0+1=1$ bit.
That of the standard basis (i.e., the root node) is clearly $2f(1/2)=2$ bits.
Therefore, the Walsh basis always wins over the standard basis.
Furthermore, the true entropy of this process is $\log \ndim=\log 2 = 1$
bit, as explained in Subsection~\ref{sec:true-entropy}.
Therefore, the mutual information of the spike process relative to the
Walsh basis is $I(\bY) = 1 - 1 = 0$ bit. 
We therefore have truly independent components for the spike process 
in this basis for $\ndim=2$, which is of course the LSDB.

\begin{hide}
The first coordinate of this vector always takes the value $\frac{1}{\sqrt{2}}$, whichever of $x_1$ or $x_2$ is selected to have a unit spike.
The second coordinate takes the values either $\frac{1}{\sqrt{2}}$ or
$\frac{-1}{\sqrt{2}}$, each with probability $\frac{1}{2}$.
We therefore have truly independent components for the spike process 
in this basis, which is of course the LSDB.

\end{hide}
\item[\underline{$\ndim=4$:}] In this case, we consider all possible
orthonormal bases in the dictionary exhaustively.
Let us mark the table of Figure~\ref{fig:haar-walsh-4} with
$+$ and $-$ signs.
\begin{figure}
\begin{center}
\begin{tabular}{|c|c|c|c|} \hline
\multicolumn{4}{|c|}{}\\ \hline
\multicolumn{2}{|c|}{$+$} &\multicolumn{2}{|c|}{$-$} \\ \hline
$++$&$-+$&$+-$&$--$  \\ \hline
\end{tabular} 
\end{center} 
\caption{The Haar-Walsh dictionary table of depth $\ndimlev=2$, i.e.,
$\ndim=4$.}
\label{fig:haar-walsh-4}
\end{figure}
We observe:
\begin{itemize}
\item each coordinate in the $-+$, $+-$, and $--$ nodes generates the same
entropy,  $g(1)=1$ bit;
\item the coordinate in the $++$ node generates $f(1)=0$ bit;
\item each coordinate in the $-$ node generates $g(1/2)=3/2$ bits;
\item each coordinate in the $+$ node generates $f(1/2)=1$ bit.
\end{itemize}
From these coordinate-wise entropy values, we can compute the entropy of
each possible basis in this dictionary as follows:
\begin{itemize}
\item the Walsh basis (the level $k=2$ basis) generates 
$1 \times 0 + 3 \times 1=3$ bits;
\item any basis using the $+$ node generates entropy larger than
$ 2 \times 1+2 \times 1=4$ bits, hence is not the LSDB;
\item any basis using the $-$ node generates entropy larger than
$2 \times \frac{3}{2}+1 \times 1=4$ bits, hence is not the LSDB;
\item the standard basis generates  $4\log 4-3\log 3 \simeq 3.24$ bits,
hence is not the LSDB.
\end{itemize}
Consequently, the LSDB for $\ndim=4$ is the Walsh basis.
This basis does not provide the truly independent components since
$I(\bY)=3-\log \ndim=3-\log 4 = 1 \neq 0$. 
\end{description}
This concludes the proof of Theorem~\ref{thm:spike-Haar-Walsh}.
{\hfill $\Box$}\vspace{1em}

\subsection{Coordinate-wise Entropy of the Spike Process}
Before proceeding to the proof of Theorems~\ref{thm:spike-OG} and \ref{thm:spike-GL}, let us consider coordinate-wise entropy of the spike process and
define some convenient quantities to characterize a basis in 
$\OG(\ndim)$ or $\GL(\ndim,\R)$.

Let us consider an invertible matrix $U=(u_{ij})_{i,j=1,\ldots,\ndim}=B^{-1} \in \GL(\ndim,\R)$, and the vector $\bY=U\bX$.
Let us consider the $i$th coordinate of $\bY$, 
$Y_{i}=\sum_{j=1}^{\ndim}u_{ij}X_{j}$.
For each realization of the spike process $\bX$, 
$Y_{i}$ takes one of the values $\{u_{ij},j=1,\ldots,\ndim \}$.
More precisely, we have $\Prob\{X_j=1\} = 1/\ndim$ and
$\Prob\{X_j=0\} = 1-1/\ndim$, for $j=1,\ldots,\ndim$.
Thus, if all $\{u_{ij},j=1, \ldots, \ndim\} $ were distinct, $Y_{i}$ would take
these values with a uniform pmf. 
But there is no particular reason that allows us to think $\{u_{ij},j=1,\ldots,\ndim\}$ are mutually distinct. 
Therefore, we shall group these values in ``classes'' of equality.
Let us introduce, for each $i\in \{1,\ldots,\ndim\}$, an integer $k(i)$
equal to the number of distinct values in the $i$th row vector
$\{u_{ij},j=1,\ldots,\ndim\}$, and the vector 
$c(i)=(\alpha_{1}(i), \ldots , \alpha_{k(i)}(i) ) \in \N^{k(i)}$, 
where each component counts the number of occurrences of each distinct value
in the $i$th row vector.
We will call $k(i)$ the \emph{class} of the $i$th row and 
$c(i)$ the \emph{index} of that row.
Clearly, we have
\bdm
1 \leq k(i) \leq \ndim \quad \mbox{and} \quad \sum_{\ell=1}^{k(i)} \alpha_{\ell}(i) = \ndim.
\edm
For example, with $\ndim=3$, if we had
\bdm
\left \{
\begin{array}{l}
Y_{1}= X_{1}+ X_{2}+X_{3} \\
Y_{2}= 5 X_{1}+2 X_{2}+2 X_{3} \\
Y_{3}= -X_{1}+X_{2}
\end{array}
\right.,
\edm
then we would get
\bdm
\left \{
\begin{array}{l}
k(1)=1, \, c(1)=(3) \\
k(2)=2, \, c(2)=(2,1) \\
k(3)=3, \, c(3)=(1,1,1)
\end{array}
\right. ,
\edm
since $\{u_{1j}\}=\{1,1,1\}$ in which we find three 1's,
$\{u_{2j}\}=\{5,2,2\}$ in which we find two 2's, one 5, and
$\{u_{3j}\}=\{-1,1,0\}$ in which we find one -1, one 1, and one 0.

Let us now examine the coordinate-wise entropy in terms of the quantities we have just defined. 
Suppose the value $u$ appears $\alpha_{\ell}(i)$ times 
in $\{u_{ij},j=1,\ldots,\ndim\}$. Then the probability of the event
$\{ Y_i=u \}$ is $\alpha_{\ell}(i)/\ndim$.
Therefore, we have
\bdm
H(Y_{i})= - \sum_{\ell=1}^{k(i)} \frac{\alpha_{\ell}(i)}{n}\ \log{\frac{\alpha_{\ell}(i)}{\ndim}}.
\edm
We shall now describe the different values that this coordinate-wise entropy 
takes as the number of distinct values and their occurrences vary.
Because the entropy is a measure of uncertainty, we can intuitively guess that
a coordinate with a small class number generates small entropy.
\begin{description}
\item[\underline{$k(i)=1$:}] This necessarily means that $c(i)=(\ndim)$, i.e.,
all the $\{u_{ij}, j=1,\ldots,\ndim\}$ are identical.
Since there is no uncertainty about this coordinate, we can expect its entropy
to be $0$.
Indeed, $H(Y_{i})=- \sum_{k=1}^{1} \frac{\ndim}{\ndim}\ \log{\frac{\ndim}{\ndim}}=0$.
\item[\underline{$k(i)=2$:}] Let us consider the link between the uncertainty
and the index $c(i)$. $k(i)=2$ means that $Y_{i}$ can take only two distinct
values. The least scattered distribution of these two values corresponds to 
the case $c(i)=(1,\ndim-1)$. 
This is also the distribution closest to the certain case $k(i)=1$ and
$c(i)=(\ndim)$.
We now show that the case $c(i)=(1,\ndim-1)$ generates the smallest entropy.
Suppose that $Y_{i}$ can take two distinct values with index $(\alpha_1,\alpha_2)$, $\alpha_1+\alpha_2=\ndim$.  In other words, $Y_i$ takes these two values
with probability $\alpha_1/\ndim$ and $\alpha_2/\ndim=1-\alpha_1/\ndim$,
respectively.  
Without loss of generality, we can assume $\alpha_1 \leq \alpha_2$. 
Therefore, the entropy of the coordinate $Y_i$ is 
\begin{eqnarray*}
H(Y_i) & = & - \left[ \frac{\alpha_1}{\ndim} \log{\frac{\alpha_1}{\ndim}}+ \frac{\alpha_2}{\ndim} \log{\frac{\alpha_2}{\ndim}} \right] \\
        & = & - \left[ \frac{\alpha_1}{\ndim} \log{\frac{\alpha_1}{\ndim}}+\(1- \frac{\alpha_1}{\ndim}\) \log{\(1-\frac{\alpha_1}{\ndim}\)} \right] \\
        & = & f\(\frac{\alpha_1}{\ndim}\),
\end{eqnarray*}
where the function $f$ is defined in \eqref{eq:f} and shown in Figure~\ref{fig:f}.
Since $\alpha_1 \leq \alpha_2$, it suffices to consider $\alpha_1$ with
$1 \leq \alpha_1 \leq \ndim/2$.  So, we have 
$1/\ndim \leq \alpha_1/\ndim \leq 1/2$, and in this interval,
$f(\alpha_1/\ndim)$ is strictly increasing.  In other words,
\bdm
f\(\frac{1}{\ndim}\) \leq f\(\frac{\alpha_1}{\ndim}\) \leq f\(\frac{1}{2}\) = 1.
\edm
Therefore, the entropy is minimal when $\alpha_1=1$ and $\alpha_2=\ndim-1$.
For $\alpha_1 \geq 2$, we clearly have $H(Y_i) \geq f(2/\ndim)$.

\item[\underline{$k(i) \geq 3$}:] To find a lower bound of
$H(Y_{i})= - \sum_{\ell=1}^{k(i)} \frac{\alpha_{\ell}(i)}{\ndim}\ \log{\frac{\alpha_{\ell}(i)}{\ndim}}$, we need the following lemma:
\begin{aismle}
\label{lem:ent-alphak}
Let $k \geq 3$ be an integer, and let $(\alpha_1,\ldots,\alpha_k)$ be
a set of strictly positive integers with
$\sum_{j=1}^{k} \alpha_j = \ndim$.
Then,
\bdm
\sum_{j=1}^{k} \frac{\alpha_j}{\ndim} \log \frac{\alpha_j}{\ndim} \leq 
-\(1+\frac{2(k-2)}{\ndim} \) f \(\frac{1}{\ndim}\).
\edm
\end{aismle}
See Appendix~A for the proof of this lemma.

Lemma~\ref{lem:ent-alphak} implies that
\bdm
H(Y_i) \geq \( 1+ \frac{2(k-2)}{\ndim} \) f\(\frac{1}{\ndim}\).
\edm
\end{description}
We can now summarize these results as the following lemma: 
\begin{aismle}
\label{lem:coordinatewise-entropy}
The coordinate-wise entropy of the spike process after transformed by
a basis in $\GL(\ndim,\R)$ can be computed or bounded as follows:
\begin{eqnarray}
\mbox{if $k(i)=1$,} & \mbox{then} & H(Y_i)=0; \label{eq:k1} \\
\mbox{if $k(i)=2$,} & \mbox{then} & H(Y_i) \left\{ \begin{array}{lll}
= & f(1/\ndim) & \mbox{ if $\, \alpha_1(i)=1$;} \\
\geq & f(2/\ndim) & \mbox{ if  $\, 2 \leq \alpha_1(i) \leq \ndim/2$;}
\end{array} \right. \label{eq:k2} \\
\mbox{if $k(i) \geq 3$,} & \mbox{then} & H(Y_i) \geq \( 1+ \frac{2(k-2)}{\ndim} \) f\(\frac{1}{\ndim}\) \geq \( 1+ \frac{2}{\ndim} \) f\(\frac{1}{\ndim}\). \label{eq:k3}
\end{eqnarray}
\end{aismle}
Let us now come back to our invertible transformation $U$; we are searching for
the LSDB among $\OG(\ndim)$ or $\GL(\ndim,\R)$.
This means that the cost of the LSDB, i.e., the sum of the coordinate-wise
entropy of the LSDB coordinates, cannot be larger than that of the standard
basis.
Therefore we will always keep the standard basis in mind as a reference 
basis with which we shall compare the performance of all other bases. 

The standard basis corresponds to $U=I_\ndim$. Every row of 
the standard basis has index $k(i)=2$ and $c(i)=(1,\ndim-1)$.
Hence the entropy cost of the standard basis is
\begin{equation}
\label{eq:ent-cost-std}
\Cost_{H}(I_\ndim \, | \, \bX) = \ndim \times f(1/\ndim)=\ndim\log \ndim - (\ndim-1)\log(\ndim-1).
\end{equation}

We saw that, assuming $k(i)>1$, $H(Y_i) \geq f(1/\ndim)$, with equality if and
only if $k(i)=2$ and $c(i)=(1,\ndim-1)$. 
Therefore a basis with $k(i)>1$ for every $i \in \{1,\ldots,\ndim\}$
has no chance to win over the standard basis, and the best thing one can do
with such a basis is to match the entropy with that of the standard basis,
i.e., a basis with $k(i)=2$ and $c(i)=(1,\ndim-1)$ for every $i$.

So, the only chance to beat the standard basis is to have some
``class 1'' rows (i.e., $k(i)=1$) in a basis.
However, we will never find an invertible matrix with more than
one class 1 rows.
Indeed, a class 1 row is necessarily proportional to 
$\bone_\ndim^T=(1,1, \ldots, 1)$,
and it is evident that more than one class 1 rows cannot exist in any 
invertible matrix.

\subsection{Proof of Theorem~\ref{thm:spike-OG}}
{\sc Proof.}
Let us start with a simple remark. 
If we assume that $B$ is an orthonormal basis, then $U=B^{-1}=B^T$. 
Hence the rows of $U$ are in fact the basis vectors of this basis.
In the case of an orthonormal matrix, the presence of one row of class 1 
imposes a constraint on the other rows, since these rows must form an 
orthonormal basis. The following lemma describes one of these constraints.
\begin{aismle}
\label{lem:SO-ones}
If $k(1)=1$, then it is impossible to have two class 2 rows with index
$(1,\ndim-1)$ in a matrix $U \in \OG(\ndim)$.
In other words,
If $k(1)=1$, then there do not exist $i_1,i_2 \in \{1,\ldots,\ndim\}$ such 
that $i_1 \neq i_2$ and $c(i_1)=c(i_2)=(1,\ndim-1)$.
\end{aismle}
The proof of this lemma can be found in Appendix~\ref{app:lem-SO-ones}.

Hence, assuming that $k(1)=1$, we can 
have at most one row of class 2 with index $(1,\ndim-1)$. 
All the other rows will be of either class $k(i)>2$ or class $k(i)=2$ with
index $(\alpha_1,\ndim-\alpha_1)$, $1 < \alpha_1 \leq \ndim/2$. 
Considering the minimization of the sum of the coordinate-wise entropy,
we must have one row of class 1 and one row of class 2 with index $(1,\ndim-1)$.  
All the other cases always increase the entropy, i.e., dependency.
From \eqref{eq:k2} and \eqref{eq:k3}, the entropy of a row with
either $k(i) > 2$ or $k(i)=2$ with index $(\alpha_1,\ndim-\alpha_1)$, 
$1 < \alpha_1 \leq \ndim/2$ is bounded from below as
\bdm
H(Y_i) \geq \min \( \( 1+\frac{2}{\ndim} \) f \(\frac{1}{\ndim}\), 
f\(\frac{2}{\ndim}\) \)
=\min \( \frac{2}{\ndim} f\(\frac{1}{\ndim}\), f\(\frac{2}{\ndim}\)-f\(\frac{1}{\ndim}\) \)+f\(\frac{1}{\ndim}\).
\edm
Therefore, combining this with \eqref{eq:k1} for $k(1)=1$ and
\eqref{eq:k2} for $\alpha_1=1$, we have
\begin{equation}
\label{eq:sum-ent-SO}
\sum_{i=1}^{\ndim} H(Y_i) \geq 0 + f \(\frac{1}{\ndim}\) 
   +(\ndim-2) \left[\min \( \frac{2}{\ndim} f\(\frac{1}{\ndim}\), f\(\frac{2}{\ndim}\)-f\(\frac{1}{\ndim}\) \) + f\(\frac{1}{\ndim}\) \right].
\end{equation}
We now use the following lemma:
\begin{aismle}
\label{lem:min-k2k3}
For $\ndim \geq 6$,
\bdm
\min \( \frac{2}{\ndim} f\(\frac{1}{\ndim}\), f\(\frac{2}{\ndim}\)-f\(\frac{1}{\ndim}\) \)=   \frac{2}{\ndim} f\(\frac{1}{\ndim}\).
\edm
\end{aismle}
{\sc Proof.}
Let us define a function: $r(x) \define x \left[ \frac{2}{x} f\(\frac{1}{x}\)-\( f\(\frac{2}{x}\)-f\(\frac{1}{x}\)\) \right]$ for $ x \geq 2 $, where
$f$ is defined in \eqref{eq:f}.
This is a continuous and monotonically-decreasing function for $x \geq 2$,
since 
\bdm
r'(x)=-\frac{2}{x^2} \log(x-1) + \log\frac{x-2}{x-1} < 0 \quad \mbox{for $x \geq 2$}.
\edm
Moreover, we have $r(5)\approx0.199$ and $r(6)\approx-0.310$, and we can 
find a zero of $r(x)$ numerically, i.e., $r(x^*)=0$ where 
$x^* \approx 5.3623$. 
These prove that this function is negative if $x \geq x^*$.
Therefore, for each integer $\ndim \geq 6$, $r(\ndim) < 0$, i.e.,
\bdm
\frac{2}{\ndim} f\(\frac{1}{\ndim}\) < f\(\frac{2}{\ndim}\)-f\(\frac{1}{\ndim}\).
\edm 
{\hfill $\Box$}

Using this lemma for $\ndim \geq 6$, \eqref{eq:sum-ent-SO} can be written as
\bdm
\sum_{i=1}^{\ndim} H(Y_i) \geq f\(\frac{1}{\ndim}\)+(\ndim-2) \left[ \frac{2}{\ndim} f\(\frac{1}{\ndim}\) + f\(\frac{1}{\ndim}\) \right]
= \left[ \frac{2(\ndim-2)}{\ndim}+\ndim-1 \right] f\(\frac{1}{\ndim}\).
\edm
Therefore, if we compare the mutual information of the new coordinates
to that of the standard basis, we have
\bdm
I(\bY)-I(\bX) \geq \left[ \frac{2(\ndim-2)}{\ndim}+\ndim-1 \right] f\(\frac{1}{\ndim}\) - \ndim f \(\frac{1}{\ndim}\) 
=\left[ \frac{2(\ndim-2)}{\ndim}-1 \right] f\(\frac{1}{\ndim}\),
\edm
That is,
\bdm
I(\bY)-I(\bX) \geq \frac{\ndim-4}{\ndim}f\(\frac{1}{\ndim}\)>0.
\edm
Thus, $B=U^{-1}=U^T$ is not the LSDB.
We have therefore proved that any orthonormal basis yields a larger mutual 
information than the standard basis for the spike process for $\ndim \geq 6$.
 
We can summarize our results so far.
\begin{itemize}
\item For $\ndim \geq 6$, the standard basis is the LSDB among $\OG(\ndim)$.
\item Any basis that yields the same mutual information as the standard basis
necessarily consists of only class 2 rows with index $(1,\ndim-1)$.
\end{itemize}

Now the question is whether there is any other basis except the standard
basis satisfying this condition.
The following lemma concludes the proof of Theorem~\ref{thm:spike-OG} for
$\ndim \geq 6$.
\begin{aismle}
\label{lem:lsdb-OG6}
For $\ndim \geq 2$, an orthonormal basis consisting of class 2 rows
with index $(1,\ndim-1)$ other than the standard basis is uniquely 
(modulo permutations and sign flips as described in Remark~\ref{rem:permute-sign}) determined as \eqref{eq:lsdb-OG5}, i.e.,
\bdm
B_{\OG(\ndim)} = \frac{1}{\ndim} \left[ \begin{array}{cccc}
 \ndim-2 &   -2   & \cdots &   -2   \\
   -2    & \ndim-2& \ddots & \vdots \\
\vdots   & \ddots & \ddots &   -2   \\
   -2    & \cdots &   -2   & \ndim-2
\end{array} \right].
\edm
\end{aismle}
The proof of this lemma can be found in Appendix~\ref{app:lem-lsdb-OG6}.
Note that this matrix becomes a permuted and sign-flipped version of
$I_2$ when $\ndim=2$, and approaches to the identity matrix
as $\ndim \to \infty$.

We now prove the particular cases, $\ndim=3,4,5$ in 
Theorem~\ref{thm:spike-OG}.
For these small values of $\ndim$, we cannot use
Lemma~\ref{lem:min-k2k3} anymore since we have
\bdm
\min \( \frac{2}{\ndim} f\(\frac{1}{\ndim}\), f\(\frac{2}{\ndim}\)-f\(\frac{1}{\ndim}\) \)= f\(\frac{2}{\ndim}\)-f\(\frac{1}{\ndim}\).
\edm
Therefore, we prove these cases by examining exhaustively all
possible indexes and the coordinate-wise entropy they generate.

\begin{description}
\item[\underline{$\ndim=3$:}]
The possible indexes are $(3)$, $(1,2)$ and $(1,1,1)$, which generates 
the following entropy values (in bits):
\begin{eqnarray*}
(3): &  H(Y_i) &=0;\\
(1,2): & H(Y_i) &=f\(\frac{1}{3}\)=-\frac{1}{3}\log\frac{1}{3}-\frac{2}{3}\log\frac{2}{3}=\log3-\frac{2}{3};\\
(1,1,1): & H(Y_i) &= 3\times \( -\frac{1}{3} \log \frac{1}{3} \) = \log 3.
\end{eqnarray*}
Once again, the only possibility for a basis to generate lower entropy than
the standard basis is to include a class 1 row with index $(3)$. 
But here we still cannot have two class 2 rows of index $(1,2)$ on top of
the class 1 row since Lemma~\ref{lem:SO-ones} still holds for $\ndim=3$. 
Therefore, the best combination is to have one row of each possible
class, which leads to the following global coordinate-wise entropy:
\bdm
0+\log 3-\frac{2}{3}+\log 3 \simeq 2.50 <  3 \log 3 - 2 \log 2 \simeq 2.75,
\edm
that is, this best possible basis is better than the standard basis.
Therefore, the LSDB is a basis including a vector of each class. 
Considering the orthonormality of the basis, we can only have
the following basis or its permuted or sign-flipped versions for $\ndim=3$:
\bdm
U^T = B = \left[
\begin{array}{ccc}
\frac{1}{\sqrt 3} & \frac{1}{\sqrt 6} & \frac{1}{\sqrt 2} \\
\frac{1}{\sqrt 3} & \frac{1}{\sqrt 6} & \frac{-1}{\sqrt 2} \\
\frac{1}{\sqrt 3} & \frac{-2}{\sqrt 6} & 0 
\end{array}
\right].
\edm

\item[\underline{$\ndim=4$:}]
The possible indexes are: $(4)$, $(1,3)$, $(2,2)$, $(1,1,2)$, and $(1,1,1,1)$,
which generate the following entropy values (in bits):
\begin{eqnarray*}
(4): & H(Y_i) &=0; \\
(1,3): & H(Y_i) &=f\(\frac{1}{4}\)=-\frac{1}{4}\log\frac{1}{4}-\frac{3}{4}\log\frac{3}{4}=2-\frac{3}{4} \log 3 \simeq 0.811;\\
(2,2): & H(Y_i) &= f\(\frac{1}{2}\) = 1;\\
(1,1,2): & H(Y_i) &= -\frac{1}{4}\log\frac{1}{4}-\frac{1}{4}\log\frac{1}{4}-\frac{1}{2}\log\frac{1}{2}=1.5;\\
(1,1,1,1): & H(Y_i) &= 4 \times \( -\frac{1}{4} \log \frac{1}{4} \)
=2.
\end{eqnarray*}
The total coordinate-wise entropy of the Walsh basis is 
$H_W \define 0+3 \times f(1/2)=3$ bits. 
We know from Theorem~\ref{thm:spike-Haar-Walsh} that $H_W$ is smaller than
that of the standard basis.
Let $U$ be an orthonormal basis, and let $\{\bb^T_i,i=1,\ldots,4 \}$ 
be its rows.
If $U$ generates smaller entropy than the Walsh basis, it necessarily
includes one class 1 row and one class 2 row with index $(1,3)$
from the same argument as the proof of Lemma~\ref{lem:SO-ones} (see Appendix~\ref{app:lem-SO-ones}).
Let us assume that $\bb^T_1$ is of class 1 and $\bb^T_2$ of class 2 with index
$(1,3)$. In other words, $c(1)=(4)$ and $c(2)=(1,3)$.
Now, $\bb^T_2$ is of the form $(a,a,a,b)$ and orthogonality with $\bb^T_1$ implies
that $ \bb^T_2$ is proportional to the vector $(1,1,1,-3)$.
Now, $U$ cannot include a class 4 row vector of index $(1,1,1,1)$.
If so, these three rows (i.e., rows of class 1, 2, and 4) generate
the entropy $0+0.811+2=2.811$ bits.
Hence, any other admissible choice for the remaining row, i.e., a class 2 row
with index $(1,3)$, which generates $0.811$ bits, or a class 2 row with
index $(2,2)$, which generates $1$ bit, or a class 3 row with index $(1,1,2)$,
which generates $1.5$ bits, ends up larger total coordinate-wise
entropy than the Walsh basis. 
Therefore we can discard these combinations immediately, and
the indexes of $\bb^T_3$ and $\bb^T_4$ must be chosen from $(2,2)$ and 
$(1,1,2)$.
If $\bb^T_3$ is of index $(2,2)$, it is of the form $(a,a,b,b)$ and 
orthogonality with $\bb^T_1$ implies that $\bb^T_3$ is proportional to 
$(a,a,-a,-a)$. 
Then, orthogonality with $\bb^T_2$ implies: $a+a-a+3a=0$, i.e., $a=0$.
Therefore the only possibility for $\bb^T_3$ and $\bb^T_4$ is to be both of index
$(1,1,2)$, each of which generates the coordinate-wise entropy $1.5$ bits.
The total coordinate-wise entropy generated by $U$ is therefore at least 
$0+0.811+2 \times 1.5 = 3.811 > 3=H_W$, hence $U^T$ is not the LSDB.
We can now conclude that the LSDB among $\OG(4)$ is the Walsh Basis.

\item[\underline{$\ndim=5$:}]
In this case, we prove that the LSDB is the standard basis or
the basis of the Householder reflection \eqref{eq:lsdb-OG5}, 
both of which consist of class 2 rows with
index $(1,4)$ only.
Indeed, using the similar argument as before, any basis generating smaller
entropy than these two bases must have a class 1 row and a class 2 row with
index $(1,4)$.  However, since the other three rows must be 
either of class 2 with different indexes or of class 3 or higher,
the total entropy of such a basis is larger than that of the standard basis
or the Householder reflection basis:
\bdm
\sum_{i=1}^5 H(Y_i) \geq 0 + f\(\frac{1}{5}\) + 3 \times f\(\frac{2}{5}\)
\simeq 3.635 > 5 \times f\(\frac{1}{5}\) \simeq 3.609.
\edm
\end{description}
This concludes the proof of Theorem~\ref{thm:spike-OG}.
{\hfill $\Box$}

\subsection{Proof of Theorem~\ref{thm:spike-GL}}
{\sc Proof.}
For the case $\Dict=\GL(\ndim,\R)$, the constraint imposed by 
Lemma~\ref{lem:SO-ones} is lifted since the rows of $U=B^{-1}$ do not have to
form an orthonormal basis anymore.  Hence we can have as many rows of class 2 
with index $(1,\ndim-1)$ as we wish, even if the first row of $U$ is of class 1.
Clearly, we still cannot have two class 1 rows because this violates the
invertibility of $U$.
Therefore, considering all these remarks and the classification of indexes 
established in the previous subsections, it is immediate to conclude that 
the combination of classes of rows leading to the smallest sum of 
coordinate-wise entropy is one row of class 1 and $\ndim-1$ rows of class 2 
with index $(1,\ndim-1)$.  This matrix reaches the lower bound for the total
coordinate-wise entropy $(\ndim-1) f(1/\ndim)$.
Considering the invertibility of the matrix with $\ndim-1$ rows of class 2,
the most general form of the admissible matrices is the following 
(modulo permutations and sign-flips mentioned in \eqref{eq:spike-LSDB-SL-analysis}):
\bdm
U_{\GL(\ndim,\R)} = B_{\GL(\ndim,\R)}^{-1} = \left[ \begin{array}{ccccccc}
a & a & \cdots & \cdots & \cdots & \cdots & a \\
b_2 & c_2 & b_2 & \cdots & \cdots & \cdots & b_2 \\
b_3 & b_3 & c_3 &  b_3  & \cdots  & \cdots & b_3 \\
\vdots & \vdots&  & \ddots &      &        & \vdots \\
\vdots & \vdots&  &        & \ddots &     & \vdots \\
b_{\ndim-1} & \cdots & \cdots & \cdots & b_{\ndim-1} & c_{\ndim-1} & b_{\ndim-1}\\
b_\ndim & \cdots & \cdots & \cdots & \cdots & b_\ndim & c_\ndim
\end{array} \right],
\edm
where $a,b_k,c_k$, $k=2,\ldots,\ndim$, must be chosen so that 
$U_{\GL(\ndim,\R)} \in \GL(\ndim,\R)$.
We can easily compute the determinant of this matrix in a similar manner
that we derived \eqref{eq:delta-a-b}:
\bdm
\det \( U_{\GL(\ndim,\R)} \)  = a \prod_{k=2}^\ndim (c_k - b_k).
\edm
Therefore, we must have $a \neq 0$ and $b_k \neq c_k$ for $k=2,\ldots, \ndim$
for $U_{\GL(\ndim,\R)}$ to be in $\GL(\ndim,\R)$.
Note that if we want to restrict the dictionary to $\SL^\pm(\ndim,\R)$,
then we must have $\det\(U_{\SL^\pm(\ndim,\R)}\)=\pm 1$, i.e.,
$a$ must satisfy $a=\pm \prod_{k=2}^\ndim (c_k - b_k)^{-1}$.

The corresponding inverse matrix \eqref{eq:spike-LSDB-SL-synthesis}
can be computed easily by elementary linear algebra, i.e., the Gauss-Jordan
method.
We show this matrix here again:
\bdm
B_{\GL(\ndim,\R)} = \left[ \begin{array}{ccccc}
\(1 + \sum_{k=2}^\ndim b_kd_k \)/a & -d_2 & -d_3 & \cdots & -d_\ndim \\
-b_2 d_2/a & d_2 & 0 & \cdots & 0 \\
-b_3 d_3/a & 0 & d_3 & \ddots & \vdots \\
\vdots & \vdots & \ddots & \ddots  & 0 \\
-b_\ndim d_\ndim/a & 0 & \cdots & 0 & d_\ndim
\end{array} \right],
\edm
where $d_k = 1/(c_k-b_k)$, $k=2,\ldots,\ndim$.
These are the LSDB pairs (analysis and synthesis respectively).
This concludes the proof of Theorem~\ref{thm:spike-GL}.
{\hfill $\Box$}

\subsection{Proof of Proposition~\ref{prop:spike-cp-OG}}
If we transform the spike process $\bX$ by the Householder reflector 
$B_{\OG(\ndim)}$ \eqref{eq:lsdb-OG5},
the number of nonzero components of $\bY=B_{\OG(\ndim)}^T\bX$ can be
easily computed as
\bdm
\Cost_0\(B_{\OG(\ndim)} \, | \, \bX \) = E \| \bY \|_0 = \ndim.
\edm
Now, let us consider the case $0 < p < 1$.
Since $\ndim \geq 2$, we have
\bdm
\Cost_p\(B_{\OG(\ndim)} \, | \, \bX \) = E \| \bY \|_p^p = \( 1 - \frac{2}{\ndim} \)^p + (\ndim-1) \( \frac{2}{\ndim} \)^p.
\edm
Let us now define the following function:
\bdm
s_p(x) \define (1-x)^p + \( \frac{2}{x}-1 \)x^p = (1-x)^p - x^p + \frac{2}{x^{1-p}},
\edm
where $0 < x = 2/\ndim \leq 1$.
Taking the derivative with respect to $x$, we have
\bdm
s'_p(x) = - p \left(\frac{1}{(1-x)^{1-p}} + \frac{1}{x^{1-p}} \right)
+ \frac{2(p-1)}{x^{2-p}} < 0,
\edm
for $0 < x < 1$ and $0 < p \leq 1$.
Therefore, in this interval, $s_p(x)$ is monotonically decreasing, and
the decisive term for the sparsity measure $\Cost_p$ is $2/x^{1-p}$.
Therefore, we have
\bdm
\lim_{\ndim \to \infty} \Cost_p\(B_{\OG(\ndim)} \, | \, \bX \) =
 \lim_{x \, \downarrow \, 0} s_p(x) = \infty \quad \mbox{for $0 < p < 1$}.
\edm
If $p=1$, then $s_1(x)=(1-x) - x + 2 = 3-2x$.
Hence, we have
\bdm
\lim_{\ndim \to \infty} \Cost_1\(B_{\OG(\ndim)} \, | \, \bX \) =
 \lim_{x \, \downarrow \, 0} s_1(x) = 3.
\edm
This completes the proof.
{\hfill $\Box$}

\subsection{Proof of Corollary~\ref{cor:spike-GL}}
\begin{figure}
\epsfxsize=0.6\textwidth
\centerline{\epsfbox{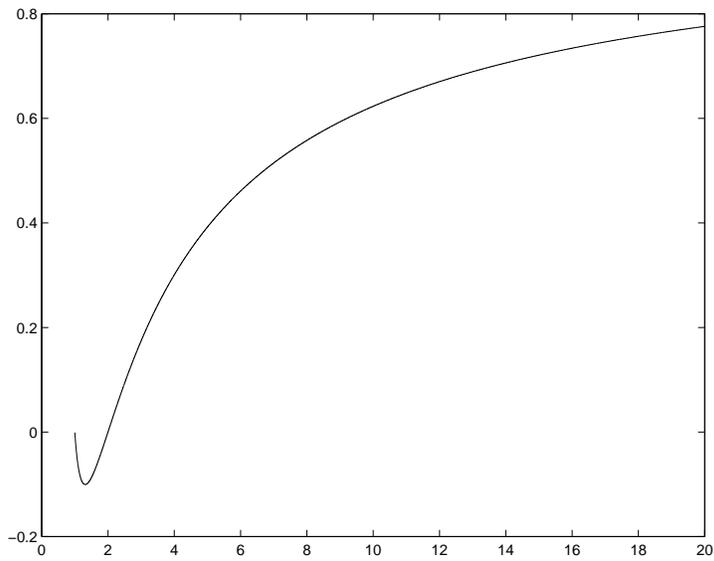}}
\caption{A plot of the function $\ln 2 \times h:x \rightarrow (x-2) \ln x -\frac{(x-1)^2}{x} \ln(x-1)$.}
\label{fig:h}
\end{figure}

\begin{figure}
\epsfxsize=0.6\textwidth
\centerline{\epsfbox{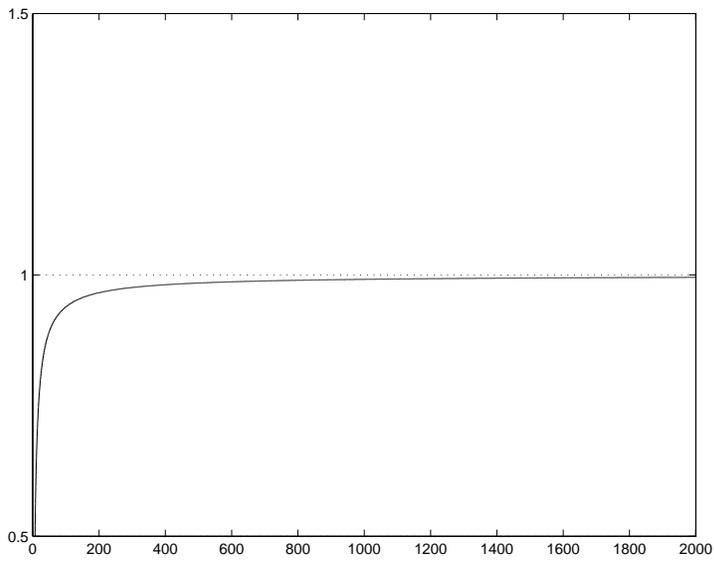}}
\caption{A plot of the function $\ln 2 \times h(x)$ for large $x$.}
\label{fig:h'}
\end{figure}

{\sc Proof.}
We now consider the mutual information of the spike process under the
LSDB pair \eqref{eq:spike-LSDB-SL-analysis} and \eqref{eq:spike-LSDB-SL-synthesis} in Theorem~\ref{thm:spike-GL}, which was proved in the previous
subsection.
Using this analysis LSDB, the mutual information of 
$\bY = B_{\GL(\ndim,\R)}^{-1} \bX$ is
\begin{eqnarray}
I(\bY) & = & -H(\bX)+\sum_{i=1}^{\ndim} H(Y_i) \nonumber \\
     & = & -\log \ndim + (\ndim-1)f \(\frac{1}{\ndim}\) \nonumber \\ 
     & = & -\log \ndim + (\ndim-1) \left[ \log \ndim -\frac{\ndim-1}{\ndim}\log(\ndim-1) \right] \nonumber \\ 
     & = & (\ndim-2)\log \ndim  -\frac{(\ndim-1)^2}{\ndim}\log(\ndim-1).
\label{eq:spike-LSDB-SL-mutual}
\end{eqnarray}
Let $h(\ndim)$ denote the last expression in \eqref{eq:spike-LSDB-SL-mutual}.
Note that $h(2)=0$, i.e., we can achieve the true independence for $\ndim=2$.
If $\ndim>2$, this function is strictly positive and monotonically
increasing as shown on Figures~\ref{fig:h} and \ref{fig:h'}.
By expanding the natural logarithm version of $h(x)$, we have
\begin{eqnarray*}
\ln 2  \times h(x)&=& (x-2) \ln x -\frac{(x-1)^2}{x} \ln(x-1) \\
    &=& (x-2) \ln x -\(x-2+\frac{1}{x}\) \(\ln x + \ln\(1- \frac{1}{x}\) \) \\
    &=& (x-2) \ln x -\(x-2+\frac{1}{x}\) \(\ln x  - \frac{1}{x}-\frac{1}{2 x^2}+o\(\frac{1}{x^2}\)\)\\
    &=& -\frac{\ln x}{x} +\(x-2+\frac{1}{x}\) \(\frac{1}{x}+\frac{1}{2 x^2}+o\(\frac{1}{x^2}\)\)\\
    &=& 1 -\frac{\ln x}{x}-\frac{3}{2x}+o\(\frac{1}{x}\) 
\end{eqnarray*}
In other words, we have established
\bdm
I(\bY) \sim \frac{1}{\ln 2}\( 1-\frac{ \ln \ndim}{\ndim} \)
\quad \mbox{as $\ndim \to \infty$.}
\edm
In other words,
\bdm
\lim_{\ndim \to \infty} I \( B_{\GL(\ndim,\R)}^{-1} \bX \) = \frac{1}{\ln 2}
=\log \e \approx 1.4427.
\edm
Therefore, for $\ndim > 2 $, there is no invertible linear transformation 
that gives truly independent coordinates for the spike process.

As for the orthonormal case, using \eqref{eq:ent-cost-std}, we have
\bdm
I\(B_{\OG(\ndim)}^T \bX \) = \ndim \log \ndim - (\ndim-1) \log (\ndim-1) - \log \ndim = (\ndim-1) \log \frac{\ndim}{\ndim-1} = \log \( 1 + \frac{1}{\ndim-1}\)^{\ndim-1}.
\edm
Now, it is easy to see
\bdm
\lim_{\ndim \to \infty} I\(B_{\OG(\ndim)}^T \bX \) = \log \e.
\edm
This completes the proof of Corollary~\ref{cor:spike-GL}.
{\hfill $\Box$}

\section{Discussion}
\label{sec:disc}
In general, sparsity and statistical independence are two completely
different concepts as an adaptive basis selection criterion, as demonstrated by
the rotations of the 2D uniform distribution in 
Section~4.
For the spike process, however, we showed that the BSB
and the LSDB \emph{can} coincide (i.e., the standard basis) if we restrict our 
basis search within $\OG(\ndim)$ with $\ndim \geq 5$. 
However, we also showed that the standard basis is not the only LSDB
in this case.  To our surprise, there exists another orthonormal 
basis \eqref{eq:lsdb-OG5} representing the Householder reflector,
which attains exactly the same level of the statistical dependence 
as the standard basis, if evaluated by the mutual information or equivalently
by the total coordinate-wise entropy $\Cost_{H}$ defined in \eqref{eq:c-lsdb}. 
Yet this LSDB does not sparsify the process at all \emph{if we measure
the sparsity by the expected $\ell^p$ norm $\Cost_p$ defined in 
\eqref{eq:cp} where $0 \leq p \leq 1$}.  It is also interesting to
note that this Householder reflector approaches to the standard basis 
as $\ndim \to \infty$.
Furthermore, if we extend our basis search to $\GL(\ndim,\R)$, then 
the LSDB and the BSB \emph{cannot} coincide.

What do these observations and the effort to prove these theorems suggest?  
First, it is clear that proving theorems on the LSDB and computing it
for more complicated stochastic processes would be much more difficult 
than the BSB.
To deal with statistical dependency, we need to consider the
``stochastics'' explicitly such as entropy and the pdf of each coordinate.
On the other hand, sparsity does not require such information.
In fact, one can even adapt the BSB for each realization rather than for
the whole realizations; see Saito et al.\ (2000, 2001) for further information
about this issue.

Second, Remark~\ref{rem:lsdb-dense} and Proposition~\ref{prop:spike-cp-OG}
cast questions on the appropriateness of
the $\ell^p$ norm ($0 \leq p < 1$) \eqref{eq:ellp-norm} as a sparsity measure.
According to this measure, the bases \eqref{eq:lsdb-OG5} and 
\eqref{eq:spike-LSDB-SL-12} provide completely ``dense'' coordinates 
for the spike process.  Yet, if we look at these basis vectors carefully,
they are very ``simple'' in the sense that at most one component differs from
all the other common components in each basis vector.  In other words,
\emph{the sparsity measured by the $\ell^p$ norm does not imply the simplicity
measured by the entropy, and vice versa.}
Therefore, if a given problem really requires the statistical independence
criterion, then we cannot replace it by the sparsity criterion in general.

Then, why the sparse basis of Olshausen and Field and the
ICA basis of Bell and Sejnowski were more or less the same?
Our interpretation to this phenomenon is the following
(see also Remark~\ref{rem:spike-convert}).
The Gabor-like functions they obtained essentially convert an input
image patch to a spike or spike-like image.  In our opinion, the image patch
size such as $16 \times 16$ pixels were crucial in their experiments.
Since those image patches are of small size,
they tend to have simpler image contents such as simple oriented edges.
It seems to us that if their algorithms were computationally feasible to accept
image patches of larger size such as $64 \times 64$ or $128 \times 128$, 
both the BSB and the LSDB would be very different from Gabor-like functions.
These large size image patches (due to rich scene variations and contents
in the patches of these sizes) cannot be converted to spikes by Gabor-like 
simple functions.

We also note that the LSDB is not guaranteed to provide the true 
statistically independent coordinates in general.
Therefore, if our interest is data compression, it seems to us that 
the pursuit of sparse representations should be encouraged rather than 
that of statistically independent representations.  This is also the view 
point indicated by Donoho (1998).
However, this does not mean to downgrade the importance of 
the statistical independence altogether.
If we want to separate mixed signals or to build empirical models of
stochastic processes for simulation purposes,
then pursuing the statistical independence should be encouraged, and 
we expect to see further interplay between these two criteria.

Finally, there are a few interesting generalizations of the spike process,
which need to be addressed in the future.
One is the spike process with varying amplitude.  The spike process
whose amplitude obeys the normal distribution was treated by
Donoho et al.\ (1998) to demonstrate the superiority of the non-Gaussian 
coding using spike location information over the Gaussian-KLB coding.
The other generalization is to randomly throw in multiple spikes to a single
realization.
If one throws in more and more spikes to one realization, the standard basis
is getting worse in terms of sparsity. It will be an interesting exercise
to consider the BSB and the LSDB for such situations.


\references

Bell, A.~J. and Sejnowski, T.~J. (1997).
The `independent components' of natural scenes are edge filters,
{\em Vision Research}, {\bf 37}, 3327--3338.

Cardoso, J.-F. (1999).
High-order contrasts for independence component analysis,
{\em Neural Computation}, {\bf 11}, 157--192.

Coifman, R.~R. and Wickerhauser, M.~V. (1992)
Entropy-based algorithms for best basis selection,
{\em IEEE Trans. Inform. Theory}, {\bf 38}, 713--719.

Cover, T.~M. and Thomas, J.~A. (1991)
{\em Elements of Information Theory}, Wiley Interscience, New York.

Daubechies, I. (1992)
{\em Ten Lectures on Wavelets}, SIAM, Philadelphia, PA.

Day, M.~M. (1940).
The spaces $L^p$ with $0 < p < 1$,
{\em Bull. Amer. Math. Soc.}, {\bf 46}, 816--823.

Donoho, D.~L. (1994).
On minimum entropy segmentation,
Technical Report, Dept. Statistics, Stanford University.

Donoho, D.~L. (1998).
Sparse components of images and optimal atomic decompositions,
Technical report, Dept. Statistics, Stanford University.

Donoho, D.~L., Vetterli, M., DeVore, R.~A., and Daubechies, I. (1998).
Data compression and harmonic analysis, {\em IEEE Trans. Inform. Theory},
{\bf 44}, 2435--2476.

Hall, P. and Morton, S.~C. (1993).
On the estimation of entropy,
{\em Ann. Inst. Statist. Math.}, {\bf 45}, 69--88.

Lin, J.-J., Saito, N., and Levine, R.~A. (2000)
An iterative nonlinear {Gaussianization} algorithm for resampling
dependent components,
{\em Proc. 2nd International Workshop on Independent Component Analysis and
Blind Signal Separation}, 245--250, IEEE.

Lin, J.-J., Saito, N., and Levine, R.~A. (2001)
An iterative nonlinear {Gaussianization} algorithm for image simulation and
synthesis, preprint.

Olshausen, B.~A. and Field, D.~J. (1996).
Emergence of simple-cell receptive field properties by learning a
sparse code for natural images,
{\em Nature}, {\bf 381}, 607--609.

Olshausen, B.~A. and Field, D.~J. (1997).
Sparse coding with an overcomplete basis set: {A} strategy employed
by {V1}?
{\em Vision Research}, {\bf 37}, 3311--3325.

Saito, N. (2000).
Local feature extraction and its applications using a library of bases,
{\em Topics in Analysis and Its Applications: Selected Theses}
(ed. R.~Coifman), 269--451, World Scientific, Singapore.

Saito, N. (2001).
Image approximation and modeling via least statistically-dependent bases,
{\em Pattern Recognition\/}, to appear.

Saito, N., Larson, B.~M., and {B\'enichou}, B., (2000).
Sparsity and statistical independence from a best-basis viewpoint,
{\em Proc. SPIE}, {\bf 4119}, {\em Wavelet Applications in Signal and
Image Processing VIII} (eds. A.~Aldroubi, A.~F.~Laine, and M.~A.~Unser),
474--486, Invited paper.

Saito, N., {B\'enichou}, B., and Larson, B.~M. (2001).
Sparsity and statistical independence in adaptive signal
representations. In preparation.

van Hateren, J.~H. and van~der Schaaf, A (1998).
Independent component filters of natural images compared with simple
cells in primary visual cortex,
{\em Proc. Royal Soc. London, Ser.~B}, {\bf 265}, 359--366.


Wickerhauser, M.~V. (1994).
{\em {Adapted Wavelet Analysis from Theory to Software}}.
A~K~Peters, Wellesley, MA.

\newpage

\appendix{Proof of Lemma~\ref{lem:ent-alphak}}
\label{app:lem-ent-alphak}
{\sc Proof.}
First we need to show another lemma as follows:

{\sc Lemma A.1.}
\label{lem:ent-p1p2}
{\it 
\hspace{0.5em} Let $p_2 \geq p_1 \geq 1$ be positive integers such that $p_1+p_2 \leq \ndim$.
Then
\bdm
\frac{p_1}{\ndim}\log \frac{p_1}{\ndim}+ \frac{p_2}{\ndim} \log \frac{p_2}{\ndim} \leq  \frac{ p_1+p_2}{\ndim} \log \frac{p_1+p_2}{\ndim} - \frac{2}{\ndim} f\(\frac{1}{\ndim}\),
\edm
where $f$ is defined in \eqref{eq:f}.
}

{\sc Proof.}
The left-hand side of the inequality can be written as
\begin{eqnarray}
\frac{p_1}{\ndim}\log \frac{p_1}{\ndim}+ \frac{p_2}{\ndim} \log \frac{p_2}{\ndim}&=&
\( \frac{p_1+p_2}{\ndim} \) \left[ \frac{p_1}{p_1+p_2} \log \frac{p_1}{\ndim}
+ \frac{p_2}{p_1+p_2} \log \frac{p_2}{\ndim} \right]  \nonumber \\
&=& \( \frac{p_1+p_2}{\ndim} \) \left[ \log \frac{p_1+p_2}{\ndim}
+  \frac{p_1}{p_1+p_2} \log \frac{p_1}{p_1+p_2}
+  \frac{p_2}{p_1+p_2} \log \frac{p_2}{p_1+p_2} \right] \nonumber  \\
&=& \( \frac{p_1+p_2}{\ndim} \) \log \( \frac{p_1+p_2}{\ndim} \) 
+   \( \frac{p_1+p_2}{\ndim} \) \left[-f \(\frac{p_1}{p_1+p_2} \) \right]  \label{alacon} 	
\end{eqnarray}
However, it is clear that
\bdm
 \frac{1}{2} \geq \frac{p_1}{p_1+p_2} \geq  \frac{1}{p_1+p_2} \geq  \frac{1}{\ndim}.
\edm
From the monotonicity of $f(x)$ for $x \in [0,1/2]$, we deduce  
\bdm
1 = f \( \frac{1}{2} \) \geq f \( \frac{p_1}{p_1+p_2} \) \geq f \( \frac{1}{\ndim} \),
\edm
which we can rewrite as 
\bdm
-1 \leq -f \( \frac{p_1}{p_1+p_2} \) \leq -f \( \frac{1}{\ndim} \).
\edm
This inequality, nonnegativity of $f$, and the assumption of this lemma
yields
\bdm
\( \frac{p_1+p_2}{\ndim} \) \left[  -f \( \frac{p_1}{p_1+p_2} \) \right] \leq -\frac{2}{\ndim} f \( \frac{1}{\ndim} \).
\edm
This inequality combined with \eqref{alacon} completes the proof of
Lemma~A.1.
{\hfill $\Box$}\vspace{1em}

Coming back to the proof of Lemma~\ref{lem:ent-alphak},
we now use induction as follows.
\begin{description}
\item[\underline{$k=3$:}]
Since $\alpha_1+\alpha_2 < \ndim$, 
we can use Lemma~A.1 to assert
\bdm
\frac{\alpha_1}{\ndim}\log \frac{\alpha_1}{\ndim}+ \frac{\alpha_2}{\ndim} \log \frac{\alpha_2}{\ndim} 
\leq  
\frac{ \alpha_1+\alpha_2}{\ndim} \log \frac{\alpha _1+\alpha_2}{\ndim} - \frac{2}{\ndim} f\(\frac{1}{\ndim}\).
\edm
Therefore,
\begin{eqnarray*}
\sum_{j=1}^3 \frac{\alpha_j}{\ndim}\log \frac{\alpha_j}{\ndim} 
&\leq & \frac{\alpha_3}{\ndim}\log \frac{\alpha_3}{\ndim} +  \frac{\alpha_1+\alpha_2}{\ndim} \log \frac{\alpha _1+\alpha_2}{\ndim} - \frac{2}{\ndim} f\(\frac{1}{\ndim}\) \\
&=& \frac{\alpha_3}{\ndim}\log \frac{\alpha_3}{\ndim}+ \(1-\frac{\alpha_3}{\ndim} \) \log \(1-\frac{\alpha_3}{\ndim} \) - \frac{2}{\ndim} f\(\frac{1}{\ndim}\)\\
&=& -f \( \frac{\alpha_3}{\ndim} \)-\frac{2}{\ndim} f \( \frac{2}{\ndim} \).
\end{eqnarray*} 
We used the fact $\sum_{j=1}^3\alpha_j=\ndim$ to derive the equality in 
the second line of the above expression. 
Since $\alpha_j \geq 1$ for $j=1,2,3$, we must have
$(\ndim-1)/\ndim > \alpha_3/\ndim \geq 1/\ndim$.
Considering the symmetry of $f(x)$ around $x=1/2$ and its behavior,
we can deduce that
\bdm
\sum_{j=1}^3 \frac{\alpha_j}{\ndim}\log \frac{\alpha_j}{\ndim} \leq -f \( \frac{1}{\ndim} \) - \frac{2}{\ndim} f \( \frac{2}{\ndim} \)
\leq - \( 1 + \frac{2}{\ndim} \) f \( \frac{1}{\ndim} \).
\edm
This nails down the case $k=3$.

\item[\underline{$k \Rightarrow k+1$}:]
Let us demonstrate that,  assuming that the formula is true for $k \geq 3$,
it is still true for $k+1$. 
We can decompose the sum 
$\sum_{j=1}^{k+1} \frac{\alpha_j}{\ndim} \log \frac{\alpha_j}{\ndim}$
in the following way:
\bdm
\sum_{j=1}^{k+1} \frac{\alpha_j}{\ndim} \log \frac{\alpha_j}{\ndim}
= \frac{\alpha_{k+1}}{\ndim} \log \frac{\alpha_{k+1}}{\ndim}
+ \frac{\alpha_k}{\ndim} \log \frac{\alpha_k}{\ndim}
+ \sum_{j=1}^{k-1} \frac{\alpha_j}{\ndim} \log \frac{\alpha_j}{\ndim}.
\edm
But once again, since $\alpha_k+\alpha_{k+1} < \ndim$, we can use 
Lemma~A.1 to reach
\begin{equation}
\label{eq:ent-alphak+1}
\frac{\alpha_{k+1}}{\ndim} \log \frac{\alpha_{k+1}}{\ndim}+ \frac{\alpha_k}{\ndim} \log \frac{\alpha_k}{\ndim} 
\leq  \frac{\alpha_{k+1}+\alpha_k}{\ndim} \log \frac{\alpha_{k+1}+\alpha_k}{\ndim} - \frac{2}{\ndim} f \( \frac{1}{\ndim} \).
\end{equation}
Let us rename a sequence $\{\alpha_j\}$ as follows:
\bdm
\beta_j = \left\{ \begin{array}{ll}
\alpha_{j+1}+\alpha_j & \mbox{if $j=k$}; \\
\alpha_j & \mbox{if $j=1,\ldots,k-1$}.
\end{array} \right.
\edm
Then, using the induction assumption, we can rewrite \eqref{eq:ent-alphak+1}
as
\bdm
\sum_{j=1}^{k+1} \frac{\alpha_j}{\ndim} \log\frac{\alpha_j}{\ndim} 
\leq \frac{\beta_k}{\ndim} \log \frac{\beta_k}{\ndim} 
 + \sum_{j=1}^{k-1} \frac{\beta_j}{\ndim} \log \frac{\beta_j}{\ndim} 
 - \frac{2}{\ndim} f \( \frac{1}{\ndim} \).
\edm
Since $\sum_{j=1}^{k} \beta_j = \sum_{j=1}^{k+1} \alpha_j = \ndim$,
we can state that
\begin{eqnarray*}
\sum_{j=1}^{k+1} \frac{\alpha_j}{\ndim} \log \frac{\alpha_j}{\ndim}
& \leq & \sum_{j=1}^{k} \frac{\beta_j}{\ndim} \log \frac{\beta_j}{\ndim}
   - \frac{2}{\ndim} f \( \frac{1}{\ndim} \)  \\
& \leq & -\( 1+\frac{2(k-2)}{\ndim} \) f \( \frac{1}{\ndim} \)
      -\frac{2}{\ndim} f \( \frac{1}{\ndim} \) \\
& = & -\( 1 + \frac{2(k-1)}{\ndim} \) f \( \frac{1}{\ndim} \).
\end{eqnarray*}
\end{description}
This concludes the proof of Lemma~\ref{lem:ent-alphak}.
{\hfill $\Box$}

\appendix{Proof of Lemma~\ref{lem:SO-ones}}
\label{app:lem-SO-ones}
{\sc Proof.}
Let us prove this lemma with \emph{reductio ad absurdum}.
Let us assume that, for example, $c(2)=c(3)=(1,\ndim-1)$.
Since the first row of $U$ is proportional to $(1,1, \ldots, 1)$, 
all the other rows must satisfy
$\sum_{j=1}^{\ndim} u_{ij}=0$ for $i=2,\ldots,\ndim$ 
because of the orthonormality condition.
Let us now consider the second row $(u_{21},\ldots,u_{2\ndim})$.
Since $c(2)=(1,\ndim-1)$, let us assume $u_{21}=a$ and 
$u_{2j}=b$, $j=2,\ldots,\ndim$ for some $a, b \in \R$.
Then the orthonormality condition implies $a+(\ndim-1)b=0$.
Since the norm of this row vector has to be one, we also have
$a^2+(\ndim-1)b^2=1$.
From these two constraints, we have
$(\ndim-1)^2 b^2+(\ndim-1) b^2 =1$.
This implies $a=\pm \sqrt{\frac{\ndim-1}{\ndim}}$ and
$b= \mp \frac{1}{\sqrt{\ndim(\ndim-1)}}$.

\noindent
As the second and third rows of $U$ must be linearly independent, 
we need to assume that the third row is $(c,d,c,\ldots,c)$ for some
$c, d \in \R$.  (We cannot assume $(d, c, \ldots, c)$ for the third row
since its inner product with the second row gives
$ad+(\ndim-1)bc=0$, which leads to $c=d$ using the values of $a$ and $b$
obtained above.)  Then, similarly to the second row, we also get
$d+(\ndim-1)c = 0$, $d^2+(\ndim-1)c^2=1$.  Thus, we have
$d=\pm a$ and $c=\pm b$.
Then, regardless of the choice of the signs for $a,b,c,d$, 
the orthogonality of the second and third rows yields
\bdm
0=(\ndim-2)b^2+2ab=(\ndim-2) \cdot \frac{1}{\ndim(\ndim-1)}-2\cdot \frac{1}{\ndim}.
\edm
This leads to
$2= \frac{\ndim-2}{\ndim-1}$, i.e., $2\ndim-2=\ndim-2$, and finally
$\ndim=0$.
This contradiction implies that the assumption made is impossible,
and proves the lemma. 
{\hfill $\Box$}\vspace{1em}

\appendix{Proof of Lemma~\ref{lem:lsdb-OG6}}
\label{app:lem-lsdb-OG6}
{\sc Proof.}
Our strategy of proving this lemma is the following.
First we will show that the LSDB selected from $\OG(\ndim)$, which consists of
only class 2 row vectors with index $(1,\ndim-1)$, must be of the form:
\begin{equation}
\label{eq:lsdb-OG6-gen}
\left[ \begin{array}{cccccc}
a_1 & b_1 & \cdots & \cdots & \cdots & b_1 \\
b_2 & a_2 & b_2 & \cdots & \cdots & b_2 \\
\vdots &  & \ddots &    & & \vdots \\
\vdots &  &  & \ddots   & & \vdots \\
b_{\ndim-1} & \cdots & \cdots & b_{\ndim-1} & a_{\ndim-1} & b_{\ndim-1}\\
b_\ndim & \cdots & \cdots & \cdots & b_\ndim & a_\ndim
\end{array} \right].
\end{equation}
where $a_k^2 + (\ndim-1)b_k^2 = 1$ for $k=1,\ldots,\ndim$.
We then derive the final form \eqref{eq:lsdb-OG5} using the orthonormality of
the row vectors of this matrix \eqref{eq:lsdb-OG6-gen}.

\noindent
Since each row is of class 2 with index $(1,\ndim-1)$, only one the entry in
a row must be different from all the other $\ndim-1$ entries. Therefore,
without loss of generality, in the $k$th row, let $a_k$ be such a 
distinguishing entry and $b_k$ be the other $\ndim-1$ entries.
Let $B=U^T$ be the LSDB under consideration.
Suppose $U$ has the $i$th and $j$th rows in which the locations of $a_i$ and
$a_j$ coincide.  
Without loss of generality (modulo row and column permutations), we can
assume that $U$ is of the following form.
\begin{equation}
\label{eq:lsdb-OG6-gen-bad}
\begin{hide}
\left[ \begin{array}{ccccc}
a_1 & b_1 & b_1 & \cdots & b_1 \\
a_2 & b_2 & b_2 & \cdots & b_2 \\
b_3 & a_3 & b_3 & \cdots & b_3 \\
\vdots &  & \ddots & \ddots & \vdots \\
b_\ndim & b_\ndim & \cdots & a_\ndim & b_\ndim
\end{array} \right].
\end{hide}
\left[ \begin{array}{cccccc}
a_1 & b_1 & \cdots & \cdots & \cdots & b_1 \\
a_2 & b_2 & \cdots & \cdots & \cdots & b_2 \\
b_3 & a_3 &   b_3  & \cdots & \cdots & b_3 \\
\vdots &  & \ddots &   & & \vdots \\
b_{\ndim-1} & \cdots & \cdots & a_{\ndim-1} & b_{\ndim-1} & b_{\ndim-1}\\
b_\ndim & \cdots & \cdots & b_\ndim & a_\ndim & b_\ndim
\end{array} \right].
\end{equation}
From the normalization condition, we must have:
\begin{equation}
\label{eq:lsdb-OG6-norm}
a_k^2 + (\ndim-1)b_k^2 = 1 \quad \mbox{for $k=1,\ldots,\ndim$}.
\end{equation}
From the orthonormality condition, $U^TU=I_\ndim$, the diagonal entries
of $U^T U$ are:
\begin{eqnarray*}
(U^T U)_{1,1} & = & 1 =  a_1^2 + a_2^2 + \sum_{j=3}^\ndim b_j^2,\\
(U^T U)_{k,k} & = & 1 =  a_k^2 + \sum_{j=1,j \neq k}^\ndim b_j^2, \quad 2 \leq k < \ndim,\\
(U^T U)_{\ndim,\ndim} & = & 1 = \sum_{j=1}^\ndim b_j^2.
\end{eqnarray*}
These imply that $a_k^2=b_k^2$ for $k \geq 3$.  
Inserting this to \eqref{eq:lsdb-OG6-norm} and noting that we must have
$a_k \neq b_k$ because of the class 2 condition, we obtain:
\begin{equation}
\label{eq:lsdb-OG6-akbk}
a_k = \pm 1/\sqrt{\ndim}, \, b_k=\mp 1/\sqrt{\ndim},  \quad \mbox{for $k \geq 3$.}
\end{equation}
Consider now the off-diagonal entry of $U^T U$, for example,
\begin{eqnarray*}
(U^T U)_{1,2} & = & 0 = a_1b_1 + a_2b_2 + a_3b_3 + b_4^2 + \cdots + b_\ndim^2,\\
(U^T U)_{1,\ndim} & = & 0 = a_1b_1 + a_2b_2 + b_3^2 + b_4^2 + \cdots + b_\ndim^2
\end{eqnarray*}
Inserting \eqref{eq:lsdb-OG6-akbk} into these, we get
\begin{eqnarray*}
a_1b_1 + a_2b_2 - \frac{1}{\ndim} + \frac{\ndim-3}{\ndim} = 0 \\
a_1b_1 + a_2b_2 + \frac{\ndim-2}{\ndim} = 0.
\end{eqnarray*}
This is a contradiction (i.e., $a_1b_1 + a_2b_2$ cannot have two different
values).  Therefore $U$ cannot have two rows where the distinguishing entries
$a_i$, $a_j$ share the same column index as \eqref{eq:lsdb-OG6-gen-bad}.
It is clear that we cannot have more than two such rows.
Therefore, $U$ must be of the form \eqref{eq:lsdb-OG6-gen}.

\noindent
Now, let us compute the entries of \eqref{eq:lsdb-OG6-gen}.
The normalization condition \eqref{eq:lsdb-OG6-norm} still holds.
Computing the diagonal entries of $U^TU=I_\ndim$, we have
\begin{equation}
\label{eq:lsdb-OG6-ortho}
(U^T U)_{k,k} = 1 = a_k^2 + \sum_{j=1,j\neq k}^\ndim b_j^2 \quad
\mbox{for $k=1,\ldots,\ndim$.}
\end{equation}
Combining \eqref{eq:lsdb-OG6-norm} and \eqref{eq:lsdb-OG6-ortho}, we have:
\bdm
\ndim b_k^2 = \sum_{j=1}^\ndim b_j^2 \quad \mbox{for $k=1,\ldots,\ndim$.}
\edm
This implies that $b_1^2=\cdots=b_\ndim^2$.
Then, from the normalization condition \eqref{eq:lsdb-OG6-norm}, we must
have $a_1^2=\cdots=a_\ndim^2$ also.
Consider now the off-diagonal entry of $U^TU$:
\bdm
(U^TU)_{1,2}=0=a_1b_1+a_2b_2+(\ndim-2)b_1^2.
\edm
Now, we must have $b_2=b_1$ or $b_2=-b_1$.  So, the above equation 
can be written as
\bdm
(U^TU)_{1,2}=0=a_1b_1 \pm a_2b_1+(\ndim-2)b_1^2.
\edm
This implies that either $b_1 = 0$ or $a_1 \pm a_2 + (\ndim-2)b_1 = 0$.
$b_1=0$ leads to $b_k=0$ and $a_k=\pm1$, i.e., the standard basis.
Let us consider now the other case, i.e., $a_1 \pm a_2 + (\ndim-2)b_1 = 0$.
Since $a_2=a_1$ or $a_2=-a_1$, these lead to either $b_1=0$ or
$2a_1 + (\ndim-2)b_1 = 0$.  The former case has been already treated.
Thus, let us proceed the latter case.
From this, we have
\begin{equation}
\label{eq:a1b1}
a_1 = \(1-\frac{\ndim}{2}\)b_1.
\end{equation}
Inserting this into \eqref{eq:lsdb-OG6-norm}, we have
\bdm
b_1^2=\frac{4}{\ndim^2}.
\edm
Consequently,
\bdm
a_1^2=1-(\ndim-1) \cdot \frac{4}{\ndim^2} = \( \frac{\ndim-2}{\ndim} \)^2.
\edm
Because of \eqref{eq:a1b1} (that is true for all $k$), we have:
\begin{equation}
a_k = \pm \frac{\ndim-2}{\ndim}, \, b_k = \mp \frac{2}{\ndim}, \quad \mbox{for $k=1,\ldots,\ndim$}.
\end{equation}
This means that the matrix $U$ must be of the following form or its permuted
and sign-flipped versions:
\bdm
\begin{hide}
B_{\OG(\ndim)} = \frac{1}{\ndim} \left[ \begin{array}{ccccc}
 \ndim-2 &   -2   & \cdots & \cdots &  -2 \\
   -2    & \ndim-2& \ddots &        & \vdots \\
\vdots   & \ddots & \ddots & \ddots & \vdots \\
\vdots   &        & \ddots & \ndim-2&  -2 \\
   -2    & \cdots & \cdots &   -2   & \ndim-2
\end{array} \right].
\end{hide}
U = \frac{1}{\ndim} \left[ \begin{array}{cccc}
 \ndim-2 &   -2   & \cdots &   -2   \\
   -2    & \ndim-2& \ddots & \vdots \\
\vdots   & \ddots & \ddots &   -2   \\
   -2    & \cdots &   -2   & \ndim-2
\end{array} \right].
\edm
It turns out that this is symmetric, so we have $B=U$.
This completes the proof of Lemma~\ref{lem:lsdb-OG6}.
{\hfill $\Box$}

\end{document}